\documentclass[aps,prx,twocolumn,superscriptaddress,nofootinbib,floatfix]{revtex4-1}

\usepackage{xcolor}
\definecolor{darkblue}{RGB}{40, 85, 120} 
\definecolor{quantumviolet}{HTML}{53257F} 
\definecolor{quantumgray}{HTML}{555555} 
\definecolor{quantumgreen}{HTML}{007474} 
\definecolor{quantumblue}{HTML}{3A5FCD} 
\definecolor{quantumpurple}{HTML}{8A2BE2} 
\definecolor{warmorange}{HTML}{FF8C42} 
\definecolor{violetpurple}{HTML}{8A2BE2} 
\definecolor{teal}{HTML}{008080} 
\definecolor{softbeige}{HTML}{ECE5D7} 
\definecolor{charcoalgray}{HTML}{2E2E2E} 
\definecolor{coralred}{HTML}{FF6F61} 
\definecolor{brightyellow}{HTML}{FFD700} 
\definecolor{cyanblue}{HTML}{1C77C3} 
\definecolor{deepbluegray}{HTML}{2C3E99} 
\definecolor{softblue}{HTML}{AFCBFF} 

\usepackage{amsmath}
\usepackage{amssymb}
\usepackage{amsthm}
\usepackage{bm}
\usepackage{braket}
\usepackage{xcolor}
\usepackage{graphicx}
\usepackage[caption=false]{subfig}
\usepackage[colorlinks, linktocpage, hypertexnames = false]{hyperref}

\hypersetup{linkcolor = {gray!80!black},
citecolor = {blue!90!black}, 
urlcolor = {gray!60!black}}

\newcommand{\orcidd}[1]{\href{https://orcid.org/#1}{\includegraphics[width=8pt]{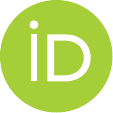}}} 

\usepackage{tikz}
\usetikzlibrary{positioning,intersections}
\usetikzlibrary{calc}
\usetikzlibrary{shapes.geometric}
\usepackage{braids}
\usepackage{tikz-cd}
\usetikzlibrary{shapes.geometric,shapes.misc}
\usetikzlibrary{arrows,matrix,calc,scopes,decorations.markings,snakes}
\usetikzlibrary{arrows.meta}
\usetikzlibrary{intersections, patterns,fit} 
\usetikzlibrary{decorations.pathreplacing,angles,quotes}

\newcommand{\black}{\color{black}}

\theoremstyle{theorem}
\newtheorem{theorem}{Theorem}

\theoremstyle{remark}

\newcommand{\Av}{\mathbb{A}}
\newcommand{\Bf}{\mathbb{B}}

\newcommand{\Hbb}{\mathbb{H}}

\newcommand{\Rep}{\mathsf{Rep}}
\newcommand{\Mod}{\mathsf{Mod}}

\newcommand{\XR}{\overset{\rightarrow}{X}}
\newcommand{\XL}{\overset{\leftarrow}{X}}

\newcommand{\ZR}{\overset{\rightarrow}{Z}}
\newcommand{\ZL}{\overset{\leftarrow}{Z}}

\newcommand{\Cocom}{\operatorname{Cocom}}

\newcommand{\cone}{(1)}
\newcommand{\ctwo}{(2)}
\newcommand{\cthree}{(3)}
\newcommand{\cfour}{(4)}

\newcommand\Tr{\operatorname{Tr}}

\newcommand{\id}{\operatorname{id}}
\newcommand{\Irr}{\operatorname{Irr}}
\newcommand{\Hom}{\operatorname{Hom}}

\newcommand{\FPdim}{\operatorname{FPdim}}

\newcommand{\Fun}{\mathsf{Fun}}
\newcommand{\Vect}{\mathsf{Vect}}

\newcommand{\cH}{\text{\usefont{OMS}{cmsy}{m}{n}H}}

\newcommand{\cZ}{\text{\usefont{OMS}{cmsy}{m}{n}Z}}

\usepackage{euscript}

\newcommand\eB           {\EuScript{B}}
\newcommand\eC           {\EuScript{C}}

\usepackage{mathptmx}

\usepackage[bbgreekl]{mathbbol} 
\DeclareMathAlphabet{\mathcal}{OMS}{cmsy}{m}{n}

\usepackage[noframe]{showframe}
\usepackage{framed}
\usepackage{lipsum}
\renewenvironment{shaded}{%
  \MakeFramed{\advance\hsize-\width \FrameRestore\FrameRestore}}%
 {\endMakeFramed}
\definecolor{shadecolor}{gray}{0.9} 

\usepackage{multirow}

\usepackage{array}
\newcolumntype{R}[1]{>{\raggedright\arraybackslash}p{#1}}

\begin{document}

\flushbottom

\title{Quantum Cluster State Model with Haagerup Fusion Category Symmetry}

\author{Zhian Jia\orcidd{0000-0001-8588-173X}}
\email{giannjia@foxmail.com}
\affiliation{Centre for Quantum Technologies, National University of Singapore, SG 117543, Singapore}
\affiliation{Department of Physics, National University of Singapore, SG 117543, Singapore}

\date{\today}

\begin{abstract}
We propose a (1+1)D lattice model, inspired by a weak Hopf algebra generalization of the cluster state model, which realizes Haagerup fusion category $\mathcal{H}_3$ symmetry and features a tensor product Hilbert space. 
The construction begins with a reconstruction of the Haagerup weak Hopf algebra \( H_3 \) from the Haagerup fusion category, ensuring that the representation category of \( H_3 \) is equivalent to \( \mathcal{H}_3 \). Utilizing the framework of symmetry topological field theory (SymTFT), we develop an ultra-thin weak Hopf quantum double model, characterized by a smooth topological boundary condition. We show that this model supports Haagerup fusion category symmetry.
Finally, we solve the ground state of the model in terms of a weak Hopf matrix product state, which serves as a natural generalization of the cluster state, embodying Haagerup fusion category symmetry.

\end{abstract}

\maketitle



\emph{Introduction.} --- 
Symmetry is one of the central topics in physics and is traditionally characterized by a group.  
From a modern perspective, symmetries are characterized by the algebraic structure of topological defects~\cite{gaiotto2015generalized,inamura2021topological,Kong2020algebraic,SchaferNameki2024ICTP,huang2023topologicalholo,bhardwaj2024lattice,freed2024topSymTFT,gaiotto2021orbifold,bhardwaj2023generalizedcharge,apruzzi2023symmetry,bhardwaj2024gappedphases,Zhang2024anomaly,Ji2020categoricalsym,chatterjee2024TopHolo,chang2019topological,bhardwaj2018finite}.  
Symmetries associated with non-invertible topological defects are known as non-invertible symmetries. For an \((n+1)\)-dimensional system, the symmetries are described by fusion \(n\)-categories~\cite{Kong2020algebraic,Bhardwaj2023non-invertible}.  
For (1+1)-dimensional systems, these symmetries are described by fusion 1-categories \(\eC\), i.e., the usual fusion categories; hence, they are also referred to as fusion category symmetries~\cite{bhardwaj2018finite,chang2019topological,thorngren2019fusion}.  

The topological defect lines (TDLs) for (1+1)D system are labeled by the simple objects in \(\eC\), and these TDLs can be fused according to the fusion rule:
\begin{equation}
    a\otimes b = \sum_{c\in \Irr(\eC)} N_{ab}^c \, c,
\end{equation}
where \(N_{ab}^c\) are non-negative integers known as fusion multiplicities. The associativity of the fusion is characterized by \(F\)-symbols. Group symmetry is a special case where the fusion category is chosen as \(\Vect_G\), the category of \(G\)-graded vector spaces \cite{etingof2016tensor}.

Non-invertible symmetries also play a crucial role in conformal field theory (CFT), particularly in the investigation of its properties and classification. In the case of rational CFTs (RCFTs), which form the foundational building blocks of general CFTs, the underlying mathematical structure is identified as modular tensor categories (MTCs) \cite{Moore1989}.
An important open question is whether a corresponding CFT can be constructed for any given MTC. It has been proposed that subfactors can be used to construct CFTs, with a conjectured correspondence between subfactors and CFTs \cite{jones1990neumann,bischoff2015relationsubfactorCFT,Bischoff2016subfactor,Xu2018subfactor,Calegari2011subfactor}. While there is significant evidence supporting this conjecture, several potential counterexamples exist, one of the most notable being the Haagerup subfactor \cite{haagerup1994principal,Asaeda1999subfactor,Asaeda2009}. There are three fusion categories arising from Haagerup subfactors, namely \( \mathcal{H}_1 \), \( \mathcal{H}_2 \), and \( \mathcal{H}_3 \), each of which is Morita equivalent to the others. 
See also Ref.~\cite{Teleman2022Haagerup} for a discussion on the special properties of the Haagerup TQFT.
Recently, there has been considerable interest in constructing lattice models that exhibit Haagerup \( \mathcal{H}_3 \) symmetry \cite{Vanhove2022Haagerup,Huang2022Haagerup,corcoran2024haagerupspin,bottini2024haagerupsymmetry,Liu2023Haagerup,Huang2022HaagerupTFT}.
The Haagerup fusion category \( \mathcal{H}_3 \) \cite{Grossman2012haagerup,Evans2011Haagerup} has a Drinfeld center \( \mathcal{Z}(\mathcal{H}_3) \) that is a modular tensor category (MTC).

\begin{figure}[b]
    \centering
    \includegraphics[width=0.9\linewidth]{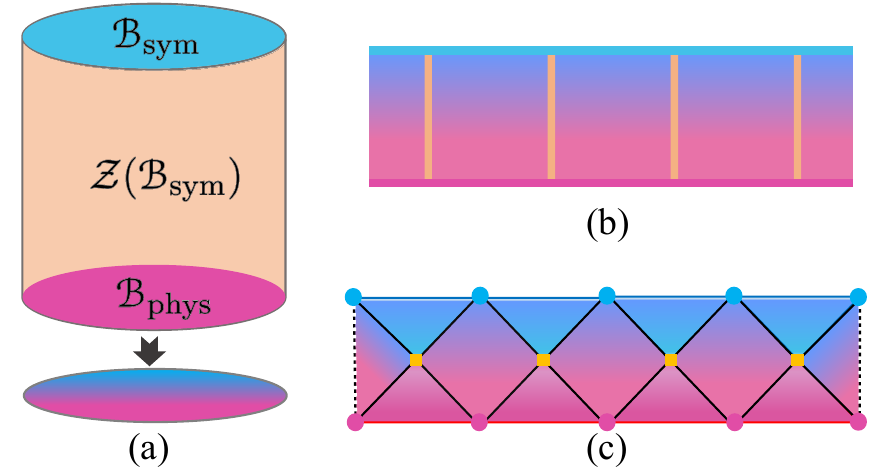}
\caption{Illustration of SymTFT sandwich and cluster ladder model. (a) The symmetry TFT consists of a symmetry boundary \( \EuScript{B}_{\rm sym} \) which encodes the fusion category symmetry; a physical boundary, which may be gapped or gapless that encodes the dynamics of the theory; the bulk is a topological field theory \( \cZ(\EuScript{B}_{\rm sym}) \). (b) Depiction of the cluster ladder model, which is an ultra-thin quantum double model with two boundaries (the qudit is put on edges), one boundary is chosen as a smooth boundary that encodes the symmetry information. If the physical boundary is chosen as a rough boundary, then the model becomes a cluster state model. (c) The chessboard representation of the cluster ladder model, where the qudits are placed on vertices, and each vertex corresponds to an edge of the quantum double model. For the cluster state model, the vertices on the physical boundary (cyan vertices) are removed, leaving only two types of local stabilizers: one for the symmetry boundary vertex operator and the other for the face operator.}
\label{fig:ladder}
\end{figure}

Given that any unitary (multi)fusion category corresponds to the representation category \( \Rep(H) \) of some weak Hopf algebra \( H \), it is a natural progression to examine non-invertible symmetries within the framework of Hopf and weak Hopf algebras~\cite{Kitaev2003,bais2003hopf,meusburger2017kitaev,Buerschaper2013a,girelli2021semidual,jia2023boundary,Jia2023weak,jia2024weakTube}. These symmetries, which are associated with structures such as Hopf, quasi-Hopf, and weak Hopf algebras, along with their module and comodule algebras, represent a substantial class of non-invertible symmetries~\cite{cordova2022snowmass,brennan2023introduction,mcgreevy2023generalized,luo2023lecture,shao2024whats,SchaferNameki2024ICTP,Bhardwaj2024lecture}.

In this work, we introduce a (1+1)D lattice model, inspired by a weak Hopf algebra generalization of the cluster state model. This model realizes the Haagerup fusion category \(\mathcal{H}_3\) symmetry and is constructed with a tensor product Hilbert space, distinguishing it significantly from models constructed using anyonic chains~\cite{Huang2022Haagerup, corcoran2024haagerupspin, bottini2024haagerupsymmetry, Liu2023Haagerup}.

The one-dimensional cluster state model is a quintessential example of an symmetry-protected topological (SPT) phase, exhibiting \(\mathbb{Z}_2 \times \mathbb{Z}_2\)  SPT order~\cite{son2012topological}. More recently, it has been proposed that the cluster state also hosts a non-invertible global symmetry, described by the fusion category \(\mathrm{Rep}(D_8)\)~\cite{seifnashri2024cluster}. 
Generalizations to finite groups are discussed in Refs.~\cite{brell2015generalized, fechisin2023noninvertible}. 
The Hopf algebraic generalization is presented in Ref.~\cite{jia2024generalized}, where it is shown that the cluster state model is, in fact, an ultra-thin quantum double model with one smooth boundary and one rough boundary. 
In a recent work~\cite{jia2024weakhopf}, we develop a general theory for the weak Hopf cluster state model, which can realize arbitrary weak Hopf symmetries (with fusion category symmetry as a special case). Building on these insights, we introduce a lattice model that exhibits Haagerup fusion category symmetry (and, more generally, Haagerup weak Hopf symmetry). The main result can be summarized as follows:

\begin{shaded}
  \begin{theorem}\label{thm:cluster}
 By applying the Tannaka-Krein reconstruction or the boundary tube algebra approach, we can recover a \(C^*\) weak Hopf algebra \(H_3\), whose representation category \(\Rep(H_3)\) is equivalent to the Haagerup fusion category \(\mathcal{H}_3\) as fusion categories. When \(H_3\) is used as input data for a weak Hopf cluster state model (or, more generally, a cluster ladder model) as described in Ref.~\cite{jia2024weakhopf}, we obtain a cluster model that exhibits Haagerup \(\mathcal{H}_3\) symmetry (see Eq.~\eqref{eq:SymHam}) and, more generally, weak Hopf symmetry \(\hat{H}_3\) (the dual weak Hopf algebra of \(H_3\)).

    The Haagerup fusion category symmetric cluster state model \(\mathbb{H}_{\rm cluster} = -\sum_{v_s} \Av_{v_s} - \sum_f \Bf_f\) possesses the following weak Hopf symmetry:
    \begin{itemize}
        \item On a closed manifold: \(\operatorname{Sym} = \operatorname{Cocom}(\hat{H}_3) \times \operatorname{Cocom}(H_3)\), where $\Cocom(H)$ is the set of all cocommutative elements in weak Hopf algebra $H$.
        \item On an open manifold: \(\operatorname{Sym} = \hat{H}_3 \times H_3\).
    \end{itemize}
    Note that \(\mathcal{H}_3 \simeq \Rep(H_3) \subset \operatorname{Cocom}(\hat{H}_3) \subset \hat{H}_3\) (here we focus on the fusion algebra of the category), implying that both cases exhibit \(\mathcal{H}_3\) symmetry.
  \end{theorem}
\end{shaded}

\vspace{0.5em}
\emph{Haagerup Fusion Category $\mathcal{H}_3$ Symmetry.} ---
The Haagerup fusion category is one of the simplest examples that does not originate from finite groups or affine Lie algebras.

The category $\mathcal{H}_3$ contains six simple objects: $\mathbf{1}$, $\alpha$, $\alpha^2$, $\rho$, ${}_{\alpha}\rho$, ${}_{\alpha^2}\rho$.
The fusion rules are summarized in Table~\ref{tab:H3fusion}, with the nontrivial rules given by:
\begin{equation}\label{eq:fusion}
    \begin{aligned}
    &  \alpha^3=\mathbf{1},\quad  {}_{\alpha}\rho =\alpha \otimes \rho,\quad  {}_{\alpha^2}\rho =\alpha^2 \otimes \rho, \\
     & \rho\otimes \alpha ={}_{\alpha^2}\rho,\quad \rho\otimes \alpha^2 ={}_{\alpha}\rho,\\
     &   \rho \otimes \rho =\mathbf{1} \oplus \rho \oplus {}_{\alpha}\rho \oplus {}_{\alpha^2}\rho.
    \end{aligned}
\end{equation}
All other fusion rules can be derived from Eq.~\eqref{eq:fusion}.
The quantum dimensions (or Frobenius-Perron dimensions, in this work, we will not distinguish them) of the simple objects are $d_{\mathbf{1}} = d_{\alpha} = d_{\alpha^2} = 1$ and $d_{\rho} = d_{{}_{\alpha}\rho} = d{{}_{\alpha^2}\rho} = \frac{3 + \sqrt{13}}{2}$. The quantum dimension of $\mathcal{H}_3$ is
\begin{equation}
    \operatorname{FPdim} \mathcal{H}_3 = \sum_{a \in \Irr (\cH_3)} d_a^2 =3(1+\frac{11+3\sqrt{13}}{2}) \approx 35.725.
\end{equation}
The $F$-symbols of $\mathcal{H}_3$ have been computed in Refs.~\cite{osborne2019fsymbolsh3fusioncategory,huang2021fsymbols}.

\begin{table}[h]
\caption{\label{tab:H3fusion}Fusion rule for the Haagerup fusion category \(\mathcal{H}_3\). In this table, the elements of the first column are fused with the elements of the first row (``column labels'' $\otimes$ ``row labels'').}
\resizebox{1.0\linewidth}{!}{%
\begin{tabular}{l|l|l|l|l|l|l}
\hline
\hline
$\otimes$  & $\mathbf{1}$ & $\alpha$ & $\alpha^2$ & $\rho$ & ${}_{\alpha}\rho$ & ${}_{\alpha^2}\rho$\\
\hline
$\mathbf{1}$  & $\mathbf{1}$ & $\alpha$ & $\alpha^2$ & $\rho$ & ${}_{\alpha}\rho$ & ${}_{\alpha^2}\rho$\\ 
\hline
$\alpha$           & $\alpha$ & $\alpha^2$ & $\mathbf{1}$ & ${}_{\alpha}\rho$ & ${}_{\alpha^2}\rho$ & $\rho$\\
\hline
$\alpha^2$         & $\alpha^2$ & $\mathbf{1}$ & $\alpha$ & ${}_{\alpha^2}\rho$ & $\rho$ & ${}_{\alpha}\rho$\\
\hline
$\rho$        & $\rho$ & ${}_{\alpha^2}\rho$ & ${}_{\alpha}\rho$ & $\mathbf{1} \oplus \rho \oplus {}_{\alpha}\rho \oplus {}_{\alpha^2}\rho$ & $\alpha^2 \oplus \rho \oplus {}_{\alpha}\rho \oplus {}_{\alpha^2}\rho$ & $\alpha \oplus \rho \oplus {}_{\alpha}\rho \oplus {}_{\alpha^2}\rho$\\
\hline
${}_{\alpha}\rho$    & ${}_{\alpha}\rho$ & $\rho$  & ${}_{\alpha^2}\rho$ & $\alpha \oplus \rho \oplus {}_{\alpha}\rho \oplus {}_{\alpha^2}\rho$ & $\mathbf{1} \oplus \rho \oplus {}_{\alpha}\rho \oplus {}_{\alpha^2}\rho$ & $\alpha^2 \oplus \rho \oplus {}_{\alpha}\rho \oplus {}_{\alpha^2}\rho$ \\
\hline
${}_{\alpha^2}\rho$& ${}_{\alpha^2}\rho$ & ${}_{\alpha}\rho$ &   $\rho$  & $\alpha^2 \oplus \rho \oplus {}_{\alpha}\rho \oplus {}_{\alpha^2}\rho$ &  $\alpha \oplus \rho \oplus {}_{\alpha}\rho \oplus {}_{\alpha^2}\rho$ & $\mathbf{1} \oplus \rho \oplus {}_{\alpha}\rho \oplus {}_{\alpha^2}\rho$ \\
\hline
\hline
\end{tabular}
}
\end{table}

The Drinfeld center $\mathcal{Z}(\mathcal{H}_3)$ is a MTC containing 12 simple objects (topological charges). These simple objects are written as $(X, \beta_{X,\bullet})$, where $X$ is a simple object in $\mathcal{H}_3$, and $\beta_{X,\bullet}$ denotes the half-braiding. Following Ref.~\cite{Hong2008exotic}, we denote these topological charges as $\mathbf{1}_{\cZ(\cH_3)},\pi_{1,2},\sigma_{1,2,3},\mu_{1,2,\cdots,6}$.
The quantum dimensions and topological spins of these topological charges are listed in Table~\ref{tab:ZH3data}. The quantum dimension of $\cZ(\cH_3)$ is 
\begin{equation}
    \FPdim(\cZ(\cH_3))= ( \operatorname{FPdim} \mathcal{H}_3 )^2.
\end{equation}
The S-matrix is given in Ref.~\cite[Section~2.2]{Hong2008exotic}.
As shown in~\cite{Teleman2022Haagerup} (also can be proved via numerical checking of Frobenius-Perron dimensions), the category $\cZ(\mathcal{H}_3)$ cannot be realized as the braided fusion category of a Chern–Simons theory associated with any compact group. This implies that it does not admit a conventional gauge theory description.

\begin{table}[h]
\caption{\label{tab:ZH3data}Simple objects of $\mathcal{Z}(\cH_3)$, specified by $(X, \beta_{X, \bullet})$ with $X \in \cH_3$ and $\beta_{X, \bullet}$ as the half-braiding, along with their quantum dimensions and topological spins.}
\resizebox{1.0\linewidth}{!}{%
\begin{tabular}{l|l|l|l|l}
\hline
\hline
$\mathcal{Z}(\cH_3)$  & $X\in \cH_3$ & $\beta_{X,\bullet}$ & quantum dimension & topological spin \\
\hline
$\mathbf{1}_{\mathcal{Z}(\cH_3)}$  & $\mathbf{1}$ & $\id$ & $1$ & $1$ \\ 
\hline
$\pi_1$           & $\pi_1=\mathbf{1} \oplus \rho \oplus {}_{\alpha}\rho \oplus {}_{\alpha^2}\rho$ & $\beta_{\pi_1,\bullet}$ & $3d_{\rho}+1$ & $1$ \\
\hline
$\pi_2$         & $\pi_2=\mathbf{1} \oplus \mathbf{1} \oplus \rho \oplus {}_{\alpha}\rho \oplus {}_{\alpha^2}\rho$ & $\beta_{\pi_2,\bullet}$ & $3d_{\rho}+2$ & $1$ \\
\hline
$\sigma_{1,2,3}$        & $\sigma=\alpha \oplus \alpha^2 \oplus \rho \oplus {}_{\alpha}\rho \oplus {}_{\alpha^2}\rho$ & $\beta^{1,2,3}_{\sigma,\bullet}$ & $3d_{\rho}+2$ & $1, e^{\pm 2\pi i/3}$ \\
\hline
$\mu_{1,2,\cdots,6}$    & $\mu= \rho \oplus {}_{\alpha}\rho \oplus {}_{\alpha^2}\rho$ & $\beta_{\mu,\bullet}^{1,2,\cdots,6}$  & $3d_{\rho}$ & $\begin{aligned}
   & e^{\pm 4\pi i/13},\\
   &e^{\pm 10\pi i/13},\\
  & e^{\pm 12\pi i/13}
\end{aligned}$ \\
\hline
\hline
\end{tabular}
}
\end{table}

\vspace{0.5em}
\emph{Reconstructing weak Hopf $H_3$ symmetry from Haagerup fusion category $\cH_3$ symmetry.} ---
Our model is based on the weak Hopf algebra $H_3$ for which $\cH_3$ is equivalent to $\Rep(H_3)$.
By definition, a complex weak Hopf algebra~\cite{BOHM1998weak} is a complex vector space $H$ equipped with algebra $(H, \mu, \eta)$ and coalgebra $(H, \Delta, \varepsilon)$ structures, along with an antipode map $S: H \to H$, satisfying specific compatibility conditions. We use Sweedler's notation for comultiplication, denoting $\Delta(h) = \sum_{(h)} h^{(1)} \otimes h^{(2)} := \sum_i h^{(1)}_i \otimes h^{(2)}_i$.
The comultiplication satisfies $\Delta(xy) = \Delta(x) \cdot \Delta(y)$, explicitly expressed as $\sum_{(xy)} (xy)^{\cone} \otimes (xy)^{\ctwo} = \sum_{(x), (y)} x^{\cone} \cdot y^{\cone} \otimes x^{\ctwo} \cdot y^{\ctwo}$. The unit satisfies weak comultiplicativity, $(\Delta \otimes \id) \circ \Delta(1_H) = (\Delta(1_H) \otimes 1_H) \cdot (1_H \otimes \Delta(1_H))$, where $\Delta(1_H) = \sum_{(1_H)} 1_H^{\cone} \otimes 1_H^{\ctwo}$. The counit satisfies weak multiplicativity, meaning that for all $x, y, z \in H$, $\varepsilon(xyz) = \sum_{(y)} \varepsilon(xy^{\cone}) \varepsilon(y^{\ctwo}z)$.
The antipode $S$ satisfies three key properties. For the left counit, $\sum_{(x)} x^{\cone} S(x^{\ctwo}) = \varepsilon_L(x) := \sum_{(1_H)} \varepsilon(1_H^{\cone} x) 1_H^{\ctwo}$. For the right counit, $\sum_{(x)} S(x^{\cone}) x^{\ctwo} = \varepsilon_R(x) := \sum_{(1_H)} 1_H^{\cone} \varepsilon(x 1_H^{\ctwo})$. Additionally, the antipode satisfies the decomposition property, $\sum_{(x)} S(x^{\cone}) x^{\ctwo} S(x^{\cthree}) = S(x)$.

The images of the left and right counit maps, \( H_L = \varepsilon_L(H) \) and \( H_R = \varepsilon_R(H) \), serve as the tensor unit in the representation category \( \Rep(H) \), which is typically a multifusion category. If \( H \) is both connected~\cite{etingof2005fusion} (also known as pure in~\cite{BOHM1998weak}) and coconnected~\cite{Nikshych2004semisimpleWHA}, then \( \Rep(H) \) is a fusion category. For more information on weak Hopf algebras and their applications to lattice models, refer to Refs.~\cite{Jia2023weak,jia2024weakhopf,jia2025quantumclusterstatespin}.

The Tannaka-Krein reconstruction offers a structured approach to obtaining a weak Hopf algebra from a fusion category. For any fusion category $\eC$, there exists a $\mathbb{C}$-algebra $A$ (which can be taken as $H_L$) such that $\eC$ can be faithfully and exactly embedded into ${_A}\Mod_A$ via a monoidal functor $F: \eC \to {_A}\Mod_A$. The algebra $H = \operatorname{End}(F)$, defined by the natural transformations of $F$, forms a weak Hopf algebra. Furthermore, the representation category $\Rep(H)$ is equivalent to $\eC$ as a unitary fusion category \cite{szlachanyi2000finite,ostrik2003module}. 
When applied to $\cH_3$, this process yields a weak Hopf algebra, denoted $H_3$, such that $\cH_3 \simeq \Rep(H_3)$.
The weak Hopf algebras whose representation categories are equivalent to $\mathcal{H}_3$ are not unique; they are Morita equivalent to one another. We will call $H_3$ Haagerup weak Hopf algebra.

An intuitive and explicit \footnote{In this work, “explicit” means that the basis and structure constants of the weak Hopf algebra can be written down concretely, and all additional data required to construct the lattice model are also explicitly known, making numerical computations feasible.}
 approach involves the concept of the boundary tube algebra \cite{Kitaev2012boundary, bridgeman2023invertible, jia2024weakTube,jia2025weakhopftubealgebra,jia2025weakhopftubealgebra}, as applied to the string-net model, where the basis elements and Haar integral can be explicitly written. For a given bulk fusion category \( \EuScript{C} \), we consider a gapped boundary described by a \( \EuScript{C} \)-module category \( {_{\EuScript{C}}}\EuScript{M} \) and construct the corresponding boundary tube algebra \( \mathbf{Tube}({_{\EuScript{C}}}\EuScript{M}) \), which is a \( C^* \) weak Hopf algebra. Setting \( \EuScript{M} = \EuScript{C} \), i.e., considering the smooth boundary, we obtain the boundary tube algebra spanned by the following basis:
\begin{equation}
\left\{  
\begin{aligned}
    \begin{tikzpicture}
        \begin{scope}
            \fill[gray!20]
                (0,1.1) arc[start angle=90, end angle=270, radius=1.1] -- 
                (0,-0.5) arc[start angle=270, end angle=90, radius=0.5] -- cycle;
        \end{scope}
           \draw[line width=0.6pt,black,->] (0,0.5)--(0,0.7);
           \draw[line width=0.6pt,black,->] (0,0.7)--(0,1.0);
           \draw[line width=0.6pt,black] (0,0.5)--(0,1.1);
        \draw[line width=.6pt,black] (0,-0.5)--(0,-1.1);
        \draw[line width=0.6pt,black,->] (0,-1.1)--(0,-0.9);
        \draw[line width=0.6pt,black,->] (0,-1.1)--(0,-0.6);
        \draw[red, line width=0.6pt] (0,0.8) arc[start angle=90, end angle=270, radius=0.8];
           \draw[red, line width=0.6pt, ->] (0,-0.8) arc[start angle=270, end angle=180, radius=0.8];
        \node[line width=0.6pt, dashed, draw opacity=0.5] at (0,1.3) {$g$};
        \node[line width=0.6pt, dashed, draw opacity=0.5] at (0,-1.3) {$c$};
        \node[line width=0.6pt, dashed, draw opacity=0.5] at (-1,0) {$a$};
        \node[line width=0.6pt, dashed, draw opacity=0.5] at (0.3,-0.7) {$\nu$};
        \node[line width=0.6pt, dashed, draw opacity=0.5] at (0,-0.2) {$e$};
        \node[line width=0.6pt, dashed, draw opacity=0.5] at (0,0.2) {$f$};
        \node[line width=0.6pt, dashed, draw opacity=0.5] at (0.3,0.7) {$\mu$};
    \end{tikzpicture}
\end{aligned}
:\quad  \begin{aligned}
    &a,c,e,f,g\in \Irr(\eC),\\
    &\mu \in \Hom(a\otimes f,g),\\
    &\nu\in \Hom(c,a\otimes e) 
\end{aligned}
\right\}.
\label{eq:tubebasis}
\end{equation}
The unit is give by
\begin{equation} \label{eq:unitH3}
1=\sum_{f,e\in \Irr(\eC)}
\begin{aligned}
    \begin{tikzpicture}
        \begin{scope}
            \fill[gray!20]
                (0,1.1) arc[start angle=90, end angle=270, radius=1.1] -- 
                (0,-0.5) arc[start angle=270, end angle=90, radius=0.5] -- cycle;
        \end{scope}
           \draw[line width=0.6pt,black,->] (0,0.5)--(0,0.9);
           \draw[line width=0.6pt,black] (0,0.5)--(0,1.1);
        \draw[line width=.6pt,black] (0,-0.5)--(0,-1.1);
        \draw[line width=0.6pt,black,->] (0,-1.1)--(0,-0.7);
        \node[line width=0.6pt, dashed, draw opacity=0.5] at (0.3,-0.7) {$e$};
        \node[line width=0.6pt, dashed, draw opacity=0.5] at (0.3,0.7) {$f$};
    \end{tikzpicture}
\end{aligned}.
\end{equation}
which means we set $a=\mathbf{1}$ and $c=e$, $f=g$ and $\mu=\id$, $\nu=\id$ in Eq.~\eqref{eq:tubebasis}.
The multiplication is defined as 
\begin{equation} \label{eq:multiConst}
    \begin{aligned}
    \begin{tikzpicture}
        \begin{scope}
            \fill[gray!20]
                (0,1.1) arc[start angle=90, end angle=270, radius=1.1] -- 
                (0,-0.5) arc[start angle=270, end angle=90, radius=0.5] -- cycle;
        \end{scope}
           \draw[line width=0.6pt,black,->] (0,0.5)--(0,0.7);
           \draw[line width=0.6pt,black,->] (0,0.7)--(0,1.0);
           \draw[line width=0.6pt,black] (0,0.5)--(0,1.1);
        \draw[line width=.6pt,black] (0,-0.5)--(0,-1.1);
        \draw[line width=0.6pt,black,->] (0,-1.1)--(0,-0.9);
        \draw[line width=0.6pt,black,->] (0,-1.1)--(0,-0.6);
        \draw[red, line width=0.6pt] (0,0.8) arc[start angle=90, end angle=270, radius=0.8];
           \draw[red, line width=0.6pt, ->] (0,-0.8) arc[start angle=270, end angle=180, radius=0.8];
        \node[line width=0.6pt, dashed, draw opacity=0.5] at (0,1.3) {$g$};
        \node[line width=0.6pt, dashed, draw opacity=0.5] at (0,-1.3) {$c$};
        \node[line width=0.6pt, dashed, draw opacity=0.5] at (-1,0) {$a$};
        \node[line width=0.6pt, dashed, draw opacity=0.5] at (0.3,-0.7) {$\nu$};
        \node[line width=0.6pt, dashed, draw opacity=0.5] at (0,-0.2) {$e$};
        \node[line width=0.6pt, dashed, draw opacity=0.5] at (0,0.2) {$f$};
        \node[line width=0.6pt, dashed, draw opacity=0.5] at (0.3,0.7) {$\mu$};
    \end{tikzpicture}
\end{aligned}
\cdot
\begin{aligned}
    \begin{tikzpicture}
        \begin{scope}
            \fill[gray!20]
                (0,1.1) arc[start angle=90, end angle=270, radius=1.1] -- 
                (0,-0.5) arc[start angle=270, end angle=90, radius=0.5] -- cycle;
        \end{scope}
           \draw[line width=0.6pt,black,->] (0,0.5)--(0,0.7);
           \draw[line width=0.6pt,black,->] (0,0.7)--(0,1.0);
           \draw[line width=0.6pt,black] (0,0.5)--(0,1.1);
        \draw[line width=.6pt,black] (0,-0.5)--(0,-1.1);
        \draw[line width=0.6pt,black,->] (0,-1.1)--(0,-0.9);
        \draw[line width=0.6pt,black,->] (0,-1.1)--(0,-0.6);
        \draw[red, line width=0.6pt] (0,0.8) arc[start angle=90, end angle=270, radius=0.8];
           \draw[red, line width=0.6pt, ->] (0,-0.8) arc[start angle=270, end angle=180, radius=0.8];
        \node[line width=0.6pt, dashed, draw opacity=0.5] at (0,1.3) {$g'$};
        \node[line width=0.6pt, dashed, draw opacity=0.5] at (0,-1.3) {$c'$};
        \node[line width=0.6pt, dashed, draw opacity=0.5] at (-1,0) {$a'$};
        \node[line width=0.6pt, dashed, draw opacity=0.5] at (0.3,-0.7) {$\nu'$};
        \node[line width=0.6pt, dashed, draw opacity=0.5] at (0,-0.2) {$e'$};
        \node[line width=0.6pt, dashed, draw opacity=0.5] at (0,0.2) {$f'$};
        \node[line width=0.6pt, dashed, draw opacity=0.5] at (0.3,0.7) {$\mu'$};
    \end{tikzpicture}
\end{aligned}=\delta_{e,c'}\delta_{f,g'}
\begin{aligned}
    \begin{tikzpicture}
        \begin{scope}
            \fill[gray!20]
                (0,1.7) arc[start angle=90, end angle=270, radius=1.7] -- 
                (0,-0.5) arc[start angle=270, end angle=90, radius=0.5] -- cycle;
        \end{scope}
           \draw[line width=0.6pt,black,->] (0,0.5)--(0,0.7);
           \draw[line width=0.6pt,black,->] (0,0.7)--(0,1.0);
                      \draw[line width=0.6pt,black,->] (0,0.5)--(0,1.6);
           \draw[line width=0.6pt,black] (0,0.5)--(0,1.7);
        \draw[line width=.6pt,black] (0,-0.5)--(0,-1.7);
        \draw[line width=0.6pt,black,->] (0,-1.1)--(0,-0.9);
        \draw[line width=0.6pt,black,->] (0,-1.1)--(0,-0.6);
        \draw[line width=0.6pt,black,->] (0,-1.7)--(0,-1.4);
        \draw[red, line width=0.6pt] (0,0.8) arc[start angle=90, end angle=270, radius=0.8];
           \draw[red, line width=0.6pt, ->] (0,-0.8) arc[start angle=270, end angle=180, radius=0.8];
        \draw[red, line width=0.6pt] (0,1.3) arc[start angle=90, end angle=270, radius=1.3];
           \draw[red, line width=0.6pt, ->] (0,-1.3) arc[start angle=270, end angle=180, radius=1.3];
                \node[line width=0.6pt, dashed, draw opacity=0.5] at (-0.3,1.3) {$\mu$};
                      \node[line width=0.6pt, dashed, draw opacity=0.5] at (-0.3,-1.3) {$\nu$};
                            \node[line width=0.6pt, dashed, draw opacity=0.5] at (0.3,1.6) {$g$};
                                      \node[line width=0.6pt, dashed, draw opacity=0.5] at (0.3,-1.6) {$c$};
        \node[line width=0.6pt, dashed, draw opacity=0.5] at (0.3,1.1) {$g'$};
        \node[line width=0.6pt, dashed, draw opacity=0.5] at (0.3,-1.1) {$c'$};
        \node[line width=0.6pt, dashed, draw opacity=0.5] at (-1,0) {$a'$};
        \node[line width=0.6pt, dashed, draw opacity=0.5] at (0.3,-0.7) {$\nu'$};
        \node[line width=0.6pt, dashed, draw opacity=0.5] at (0,-0.2) {$e'$};
        \node[line width=0.6pt, dashed, draw opacity=0.5] at (0,0.2) {$f'$};
        \node[line width=0.6pt, dashed, draw opacity=0.5] at (0.3,0.7) {$\mu'$};
                \node[line width=0.6pt, dashed, draw opacity=0.5] at (-1.5,0) {$a$};
    \end{tikzpicture}
\end{aligned}.
\end{equation}
After evaluating the diagram on the right-hand side using \( F \)-moves, parallel moves, and loop moves \cite{jia2024weakTube}, we obtain a linear combination of the basis elements in Eq.~\eqref{eq:tubebasis}.
The counit is given by
\begin{equation}
    \varepsilon\left( 
    \begin{aligned}
    \begin{tikzpicture}
        \begin{scope}
            \fill[gray!20]
                (0,1.1) arc[start angle=90, end angle=270, radius=1.1] -- 
                (0,-0.5) arc[start angle=270, end angle=90, radius=0.5] -- cycle;
        \end{scope}
           \draw[line width=0.6pt,black,->] (0,0.5)--(0,0.7);
           \draw[line width=0.6pt,black,->] (0,0.7)--(0,1.0);
           \draw[line width=0.6pt,black] (0,0.5)--(0,1.1);
        \draw[line width=.6pt,black] (0,-0.5)--(0,-1.1);
        \draw[line width=0.6pt,black,->] (0,-1.1)--(0,-0.9);
        \draw[line width=0.6pt,black,->] (0,-1.1)--(0,-0.6);
        \draw[red, line width=0.6pt] (0,0.8) arc[start angle=90, end angle=270, radius=0.8];
           \draw[red, line width=0.6pt, ->] (0,-0.8) arc[start angle=270, end angle=180, radius=0.8];
        \node[line width=0.6pt, dashed, draw opacity=0.5] at (0,1.3) {$g$};
        \node[line width=0.6pt, dashed, draw opacity=0.5] at (0,-1.3) {$c$};
        \node[line width=0.6pt, dashed, draw opacity=0.5] at (-1,0) {$a$};
        \node[line width=0.6pt, dashed, draw opacity=0.5] at (0.3,-0.7) {$\nu$};
        \node[line width=0.6pt, dashed, draw opacity=0.5] at (0,-0.2) {$e$};
        \node[line width=0.6pt, dashed, draw opacity=0.5] at (0,0.2) {$f$};
        \node[line width=0.6pt, dashed, draw opacity=0.5] at (0.3,0.7) {$\mu$};
    \end{tikzpicture}
\end{aligned}
    \right) = 
    \frac{
    \delta_{f,e} \delta_{c,g} }{d_g}
    \begin{aligned}
    \begin{tikzpicture}
        \begin{scope}
            \fill[gray!20]
                (0,1.1) arc[start angle=90, end angle=270, radius=1.1] -- 
                (0,-0.5) arc[start angle=270, end angle=90, radius=0.5] -- cycle;
        \end{scope}
           \draw[line width=0.6pt,black,->] (0,0.7)--(0,1.0);
           \draw[line width=0.6pt,black] (0,0.5)--(0,1.1);
        \draw[line width=.6pt,black] (0,-1.1)--(0,1.1);
        \draw[line width=0.6pt,black,->] (0,-1.1)--(0,-0.9);
        \draw[line width=0.6pt,black,->] (0,-1.1)--(0,0.2);
        \draw[red, line width=0.6pt] (0,0.8) arc[start angle=90, end angle=270, radius=0.8];
           \draw[red, line width=0.6pt, ->] (0,-0.8) arc[start angle=270, end angle=180, radius=0.8];
 \draw[black, line width=0.6pt] (0,-1.1) arc[start angle=-90, end angle=90, radius=1.1];
        \node[line width=0.6pt, dashed, draw opacity=0.5] at (0,1.3) {$g$};
        \node[line width=0.6pt, dashed, draw opacity=0.5] at (-1,0) {$a$};
        \node[line width=0.6pt, dashed, draw opacity=0.5] at (0.3,-0.7) {$\nu$};
        \node[line width=0.6pt, dashed, draw opacity=0.5] at (0.3,-0.2) {$e$};
        \node[line width=0.6pt, dashed, draw opacity=0.5] at (0.3,0.7) {$\mu$};
    \end{tikzpicture}
\end{aligned}.
\end{equation}
After performing the topological evaluation of the right-hand side, the result will be a real number $\delta_{f,e}\delta_{c,g} \delta_{\mu \nu}\sqrt{\frac{d_a d_e}{d_g}}$.
The comultiplication is defined as
\begin{equation}
      \Delta\left( 
    \begin{aligned}
    \begin{tikzpicture}
        \begin{scope}
            \fill[gray!20]
                (0,1.1) arc[start angle=90, end angle=270, radius=1.1] -- 
                (0,-0.5) arc[start angle=270, end angle=90, radius=0.5] -- cycle;
        \end{scope}
           \draw[line width=0.6pt,black,->] (0,0.5)--(0,0.7);
           \draw[line width=0.6pt,black,->] (0,0.7)--(0,1.0);
           \draw[line width=0.6pt,black] (0,0.5)--(0,1.1);
        \draw[line width=.6pt,black] (0,-0.5)--(0,-1.1);
        \draw[line width=0.6pt,black,->] (0,-1.1)--(0,-0.9);
        \draw[line width=0.6pt,black,->] (0,-1.1)--(0,-0.6);
        \draw[red, line width=0.6pt] (0,0.8) arc[start angle=90, end angle=270, radius=0.8];
           \draw[red, line width=0.6pt, ->] (0,-0.8) arc[start angle=270, end angle=180, radius=0.8];
        \node[line width=0.6pt, dashed, draw opacity=0.5] at (0,1.3) {$g$};
        \node[line width=0.6pt, dashed, draw opacity=0.5] at (0,-1.3) {$c$};
        \node[line width=0.6pt, dashed, draw opacity=0.5] at (-1,0) {$a$};
        \node[line width=0.6pt, dashed, draw opacity=0.5] at (0.3,-0.7) {$\nu$};
        \node[line width=0.6pt, dashed, draw opacity=0.5] at (0,-0.2) {$e$};
        \node[line width=0.6pt, dashed, draw opacity=0.5] at (0,0.2) {$f$};
        \node[line width=0.6pt, dashed, draw opacity=0.5] at (0.3,0.7) {$\mu$};
    \end{tikzpicture}
\end{aligned}
    \right)
    =
    \sum_{k, l,\zeta} \sqrt{ \frac{d_l}{d_kd_a} }   \begin{aligned}
    \begin{tikzpicture}
        \begin{scope}
            \fill[gray!20]
                (0,1.1) arc[start angle=90, end angle=270, radius=1.1] -- 
                (0,-0.5) arc[start angle=270, end angle=90, radius=0.5] -- cycle;
        \end{scope}
           \draw[line width=0.6pt,black,->] (0,0.5)--(0,0.7);
           \draw[line width=0.6pt,black,->] (0,0.7)--(0,1.0);
           \draw[line width=0.6pt,black] (0,0.5)--(0,1.1);
        \draw[line width=.6pt,black] (0,-0.5)--(0,-1.1);
        \draw[line width=0.6pt,black,->] (0,-1.1)--(0,-0.9);
        \draw[line width=0.6pt,black,->] (0,-1.1)--(0,-0.6);
        \draw[red, line width=0.6pt] (0,0.8) arc[start angle=90, end angle=270, radius=0.8];
           \draw[red, line width=0.6pt, ->] (0,-0.8) arc[start angle=270, end angle=180, radius=0.8];
        \node[line width=0.6pt, dashed, draw opacity=0.5] at (0,1.3) {$g$};
        \node[line width=0.6pt, dashed, draw opacity=0.5] at (0,-1.3) {$l$};
        \node[line width=0.6pt, dashed, draw opacity=0.5] at (-1,0) {$a$};
        \node[line width=0.6pt, dashed, draw opacity=0.5] at (0.3,-0.7) {$\zeta$};
        \node[line width=0.6pt, dashed, draw opacity=0.5] at (0,-0.2) {$k$};
        \node[line width=0.6pt, dashed, draw opacity=0.5] at (0,0.2) {$f$};
        \node[line width=0.6pt, dashed, draw opacity=0.5] at (0.3,0.7) {$\mu$};
    \end{tikzpicture}
\end{aligned} \otimes     
\begin{aligned}
    \begin{tikzpicture}
        \begin{scope}
            \fill[gray!20]
                (0,1.1) arc[start angle=90, end angle=270, radius=1.1] -- 
                (0,-0.5) arc[start angle=270, end angle=90, radius=0.5] -- cycle;
        \end{scope}
           \draw[line width=0.6pt,black,->] (0,0.5)--(0,0.7);
           \draw[line width=0.6pt,black,->] (0,0.7)--(0,1.0);
           \draw[line width=0.6pt,black] (0,0.5)--(0,1.1);
        \draw[line width=.6pt,black] (0,-0.5)--(0,-1.1);
        \draw[line width=0.6pt,black,->] (0,-1.1)--(0,-0.9);
        \draw[line width=0.6pt,black,->] (0,-1.1)--(0,-0.6);
        \draw[red, line width=0.6pt] (0,0.8) arc[start angle=90, end angle=270, radius=0.8];
           \draw[red, line width=0.6pt, ->] (0,-0.8) arc[start angle=270, end angle=180, radius=0.8];
        \node[line width=0.6pt, dashed, draw opacity=0.5] at (0,1.3) {$l$};
        \node[line width=0.6pt, dashed, draw opacity=0.5] at (0,-1.3) {$c$};
        \node[line width=0.6pt, dashed, draw opacity=0.5] at (-1,0) {$a$};
        \node[line width=0.6pt, dashed, draw opacity=0.5] at (0.3,-0.7) {$\nu$};
        \node[line width=0.6pt, dashed, draw opacity=0.5] at (0,-0.2) {$e$};
        \node[line width=0.6pt, dashed, draw opacity=0.5] at (0,0.2) {$k$};
        \node[line width=0.6pt, dashed, draw opacity=0.5] at (0.3,0.7) {$\zeta$};
    \end{tikzpicture}
\end{aligned}.
\end{equation}
The antipode map is defined as
\begin{equation} \label{eq:antipodeTH3}
         S\left( 
    \begin{aligned}
    \begin{tikzpicture}
        \begin{scope}
            \fill[gray!20]
                (0,1.1) arc[start angle=90, end angle=270, radius=1.1] -- 
                (0,-0.5) arc[start angle=270, end angle=90, radius=0.5] -- cycle;
        \end{scope}
           \draw[line width=0.6pt,black,->] (0,0.5)--(0,0.7);
           \draw[line width=0.6pt,black,->] (0,0.7)--(0,1.0);
           \draw[line width=0.6pt,black] (0,0.5)--(0,1.1);
        \draw[line width=.6pt,black] (0,-0.5)--(0,-1.1);
        \draw[line width=0.6pt,black,->] (0,-1.1)--(0,-0.9);
        \draw[line width=0.6pt,black,->] (0,-1.1)--(0,-0.6);
        \draw[red, line width=0.6pt] (0,0.8) arc[start angle=90, end angle=270, radius=0.8];
           \draw[red, line width=0.6pt, ->] (0,-0.8) arc[start angle=270, end angle=180, radius=0.8];
        \node[line width=0.6pt, dashed, draw opacity=0.5] at (0,1.3) {$g$};
        \node[line width=0.6pt, dashed, draw opacity=0.5] at (0,-1.3) {$c$};
        \node[line width=0.6pt, dashed, draw opacity=0.5] at (-1,0) {$a$};
        \node[line width=0.6pt, dashed, draw opacity=0.5] at (0.3,-0.7) {$\nu$};
        \node[line width=0.6pt, dashed, draw opacity=0.5] at (0,-0.2) {$e$};
        \node[line width=0.6pt, dashed, draw opacity=0.5] at (0,0.2) {$f$};
        \node[line width=0.6pt, dashed, draw opacity=0.5] at (0.3,0.7) {$\mu$};
    \end{tikzpicture}
\end{aligned}
    \right)
    =   \frac{d_f}{d_g} \begin{aligned}
    \begin{tikzpicture}
        \begin{scope}
            \fill[gray!20]
                (0,1.1) arc[start angle=90, end angle=270, radius=1.1] -- 
                (0,-0.5) arc[start angle=270, end angle=90, radius=0.5] -- cycle;
        \end{scope}
           \draw[line width=0.6pt,black,->] (0,0.5)--(0,0.7);
           \draw[line width=0.6pt,black,->] (0,0.7)--(0,1.0);
           \draw[line width=0.6pt,black] (0,0.5)--(0,1.1);
        \draw[line width=.6pt,black] (0,-0.5)--(0,-1.1);
        \draw[line width=0.6pt,black,->] (0,-1.1)--(0,-0.9);
        \draw[line width=0.6pt,black,->] (0,-1.1)--(0,-0.6);
        \draw[red, line width=0.6pt] (0,0.8) arc[start angle=90, end angle=270, radius=0.8];
           \draw[red, line width=0.6pt, ->] (0,-0.8) arc[start angle=270, end angle=180, radius=0.8];
        \node[line width=0.6pt, dashed, draw opacity=0.5] at (0,1.3) {$e$};
        \node[line width=0.6pt, dashed, draw opacity=0.5] at (0,-1.3) {$f$};
        \node[line width=0.6pt, dashed, draw opacity=0.5] at (-1,0) {$\bar{a}$};
        \node[line width=0.6pt, dashed, draw opacity=0.5] at (0.3,-0.7) {$\mu$};
        \node[line width=0.6pt, dashed, draw opacity=0.5] at (0,-0.2) {$g$};
        \node[line width=0.6pt, dashed, draw opacity=0.5] at (0,0.2) {$c$};
        \node[line width=0.6pt, dashed, draw opacity=0.5] at (0.3,0.7) {$\nu$};
    \end{tikzpicture}
\end{aligned}.
\end{equation}
It can be proved that this is a \( C^* \) weak Hopf algebra \cite{jia2024weakTube,jia2025weakhopftubealgebra}.
The Haar integral is of the form:
\begin{equation}\label{eq:HaarTube}
    \lambda = \frac{1}{\operatorname{rank} \eC} \sum_{a,x,y,\mu} 
    \sqrt{\frac{d_a}{d_x^3 d_y}}
        \begin{aligned}
    \begin{tikzpicture}
        \begin{scope}
            \fill[gray!20]
                (0,1.1) arc[start angle=90, end angle=270, radius=1.1] -- 
                (0,-0.5) arc[start angle=270, end angle=90, radius=0.5] -- cycle;
        \end{scope}
           \draw[line width=0.6pt,black,->] (0,0.5)--(0,0.7);
           \draw[line width=0.6pt,black,->] (0,0.7)--(0,1.0);
           \draw[line width=0.6pt,black] (0,0.5)--(0,1.1);
        \draw[line width=.6pt,black] (0,-0.5)--(0,-1.1);
        \draw[line width=0.6pt,black,->] (0,-1.1)--(0,-0.9);
        \draw[line width=0.6pt,black,->] (0,-1.1)--(0,-0.6);
        \draw[red, line width=0.6pt] (0,0.8) arc[start angle=90, end angle=270, radius=0.8];
           \draw[red, line width=0.6pt, ->] (0,-0.8) arc[start angle=270, end angle=180, radius=0.8];
        \node[line width=0.6pt, dashed, draw opacity=0.5] at (0,1.3) {$y$};
        \node[line width=0.6pt, dashed, draw opacity=0.5] at (0,-1.3) {$y$};
        \node[line width=0.6pt, dashed, draw opacity=0.5] at (-1,0) {$a$};
        \node[line width=0.6pt, dashed, draw opacity=0.5] at (0.3,-0.7) {$\mu$};
        \node[line width=0.6pt, dashed, draw opacity=0.5] at (0,-0.3) {$x$};
        \node[line width=0.6pt, dashed, draw opacity=0.5] at (0,0.3) {$x$};
        \node[line width=0.6pt, dashed, draw opacity=0.5] at (0.3,0.7) {$\mu$};
    \end{tikzpicture}
\end{aligned}.
\end{equation}
A detailed is provided in Ref.~\cite{jia2025weakhopftubealgebra}. Since the dual weak Hopf algebra can be regarded as a \(\mathbf{Tube}(\eC_{\eC})\), the Haar measure \(\Lambda\) (Haar integral of the dual weak Hopf algebra) can be constructed similarly \cite{jia2025weakhopftubealgebra} (see also Ref.~\cite{bridgeman2023invertible} for a slightly different construction, where the structure constants of the boundary tube algebra are chosen differently, although the two algebras remain Morita equivalent).

The fusion category \( \EuScript{C} = \Fun_{\EuScript{C}}(\EuScript{C}, \EuScript{C}) \), interpreted as the boundary phase, can be embedded into the representation category of the tube algebra, \( \Rep(\mathbf{Tube}({_{\EuScript{C}}}\EuScript{C})) \). It is widely believed that \( \EuScript{C} \simeq \Rep(\mathbf{Tube}({_{\EuScript{C}}}\EuScript{C})) \) \cite{Kitaev2012boundary, bridgeman2023invertible}. In Ref.~\cite{jia2025weakhopftubealgebra}, we will argue that this is, in fact, an embedding \( \EuScript{C} \hookrightarrow \Rep(\mathbf{Tube}({_{\EuScript{C}}}\EuScript{C})) \) in the general case. While the full equivalence still requires further investigation in the string diagrammatic setting, a recent proof using internal Homs has been claimed in Ref.~\cite{bai2025weakhopf}.
Consequently, the fusion category \( \EuScript{C}_{\EuScript{M}}^{\vee} := \Fun_{\EuScript{C}}(\EuScript{M}, \EuScript{M}) \) can be realized via the weak Hopf algebra \( \mathbf{Tube}({_{\EuScript{C}}}\EuScript{C}) \).

For further details about the boundary tube algebra, we refer the reader to Refs.~\cite{jia2024weakTube, jia2025weakhopftubealgebra}. When applied to \( \mathcal{H}_3 \), we denote the corresponding tube algebra as \( H_3 = \mathbf{Tube}({_{\mathcal{H}_3}}\mathcal{H}_3) \), with a slight abuse of notation.
Notice that $\cH_3$ is a multiplicity-free fusion category, meaning that $\operatorname{dim}\operatorname{Hom}(a \otimes b, c) = 0$ or $1$, so the vertex labels can be omitted.


An interesting property of Haagerup boundary tube algebra $H_3$ is that its dual can also be understood from the perspective of boundary tube algebras~\cite{jia2024weakTube,jia2025weakhopftubealgebra}. Recall that dual  algebra $\hat{H}_{3}$ consists of functionals over $H_3$:
\begin{equation}
  \hat{H}_{3}=\{\varphi:   H_3 \to \mathbb{C}\}.
\end{equation}
This is also a weak Hopf algebra with the structure given by \emph{canonical pairing} $\langle \psi, h\rangle:=\psi(h)$, its weak Hopf algebra structure is defined as
\begin{align}
&	\langle \hat{\mu}(\varphi\otimes \psi),x\rangle=\langle\varphi\otimes \psi, \Delta(x)\rangle,\label{eq:MultiplicationBarA}\\
&	\langle \hat{\eta} (1),x\rangle= \varepsilon(x),   \; \text{i.e.,}\; \hat{1}=\varepsilon,\\
&   \langle \hat{\Delta}(\varphi), x\otimes y \rangle=\langle \varphi, \mu(x\otimes y)\rangle,\label{eq:ComultiplicationBarA} \\
&    \hat{\varepsilon}(\varphi)=\langle \varphi, \eta(1)\rangle,\\
&    \langle \hat{S}(\varphi),x\rangle =\langle \varphi, S(x)\rangle,\label{eq:AntipodeBarA}
\end{align}
where $\hat{\mu}$, $\hat{\eta}$, $\hat{\Delta}$, $\hat{\varepsilon}$, $\hat{S}$ are structure maps for $\hat{H}_3$.

The category $\cH_3$ can also be regarded as a right $\cH_3$-module category, we may likewise consider the corresponding construction from the right module structure and obtain $\mathbf{Tube}({\cH_3}_{\cH_3})$ with bases:
\begin{equation}
\left\{  
\begin{aligned}
    \begin{tikzpicture}
        \begin{scope}
            \fill[gray!20]
                (0,1.1) arc[start angle=90, end angle=-90, radius=1.1] -- 
                (0,-0.5) arc[start angle=-90, end angle=90, radius=0.5] -- cycle;
        \end{scope}
          \draw[line width=0.6pt,black,->] (0,0.5)--(0,0.7);
           \draw[line width=0.6pt,black,->] (0,0.7)--(0,1.0);
           \draw[line width=0.6pt,black] (0,0.5)--(0,1.1);
        \draw[line width=.6pt,black] (0,-0.5)--(0,-1.1);
        \draw[line width=0.6pt,black,->] (0,-1.1)--(0,-0.9);
        \draw[line width=0.6pt,black,->] (0,-1.1)--(0,-0.6);
        \draw[blue, line width=0.6pt] (0,-0.8) arc[start angle=-90, end angle=90, radius=0.8];
          \draw[blue, line width=0.6pt, ->] (0,-0.8) arc[start angle=-90, end angle=0, radius=0.8];
        \node[line width=0.6pt, dashed, draw opacity=0.5] at (0,1.3) {$g$};
        \node[line width=0.6pt, dashed, draw opacity=0.5] at (0,-1.3) {$c$};
        \node[line width=0.6pt, dashed, draw opacity=0.5] at (1,0) {$b$};
        \node[line width=0.6pt, dashed, draw opacity=0.5] at (-0.3,-0.7) {$\nu$};
        \node[line width=0.6pt, dashed, draw opacity=0.5] at (0,-0.2) {$e$};
        \node[line width=0.6pt, dashed, draw opacity=0.5] at (0,0.2) {$f$};
        \node[line width=0.6pt, dashed, draw opacity=0.5] at (-0.3,0.7) {$\mu$};
    \end{tikzpicture}
\end{aligned}
:\,  \begin{aligned}
    &b,c,e,f,g\in \Irr(\cH_3),\\
    &\mu \in \Hom_{\cH_3}(f\otimes b,g),\\
    &\nu\in \Hom_{\cH_3}(c,e\otimes b) 
\end{aligned}
\right\}.
\label{eq:tubebasis}
\end{equation}
As proved in~\cite{jia2025weakhopftubealgebra} (see Lemma~14 and Proposition~17 therein), $\mathbf{Tube}({\mathcal{H}_3}_{\mathcal{H}_3})$ can be embedded into the dual algebra $\hat{H}_3$. Since these algebras have the same dimension, they are in fact isomorphic: 
$
\mathbf{Tube}({\mathcal{H}_3}_{\mathcal{H}_3})^{\rm cop} \cong \hat{H}_3.
$
here by ``com'' we mean the comultiplication is opposite, $\Delta^{\rm op}(x):=\sum_{(x)} x^{(2)}\otimes x^{(1)}$.
The \emph{skew-pairing}  between $\mathbf{Tube}({\mathcal{H}_3}_{\mathcal{H}_3})$ and $H_3$ is given by (we omit arrows here to avoid cluttering the equation):
\begin{equation} \label{eq:pairing_def}
\begin{aligned}
  &  p\left(\begin{aligned}
    \begin{tikzpicture}[scale=0.65]
        \begin{scope}
            \fill[gray!15]
                (0,-1.5) arc[start angle=-90, end angle=90, radius=1.5] -- 
                (0,0.5) arc[start angle=90, end angle=-90, radius=0.5] -- cycle;
        \end{scope}
        \draw[line width=0.6pt,black] (0,0.5)--(0,1.5);
        \draw[line width=.6pt,black] (0,-0.5)--(0,-1.5);
        \draw[blue, line width=0.6pt] (0,1.1) arc[start angle=90, end angle=-90, radius=1.1];
        \node[ line width=0.6pt, dashed, draw opacity=0.5] (a) at (0,-1.7){$\scriptstyle s$};
        \node[ line width=0.6pt, dashed, draw opacity=0.5] (a) at (1.3,0){$\scriptstyle b$};
        \node[ line width=0.6pt, dashed, draw opacity=0.5] (a) at (0.25,-0.7){$\scriptstyle t$};
        \node[ line width=0.6pt, dashed, draw opacity=0.5] (a) at (-0.25,-1.15){$\scriptstyle \mu$};
        \node[ line width=0.6pt, dashed, draw opacity=0.5] (a) at (0.25,0.7){$\scriptstyle u$};
        \node[ line width=0.6pt, dashed, draw opacity=0.5] (a) at (0,1.7){$\scriptstyle v$};
        \node[ line width=0.6pt, dashed, draw opacity=0.5] (a) at (-0.25,1.15){$\scriptstyle \gamma$};
    \end{tikzpicture}
\end{aligned}\;,\;
\begin{aligned}
    \begin{tikzpicture}[scale=0.65]
        \begin{scope}
            \fill[gray!15]
                (0,1.5) arc[start angle=90, end angle=270, radius=1.5] -- 
                (0,-0.5) arc[start angle=270, end angle=90, radius=0.5] -- cycle;
        \end{scope}
        \draw[line width=0.6pt,black] (0,0.5)--(0,1.5);
        \draw[line width=.6pt,black] (0,-0.5)--(0,-1.5);
        \draw[red, line width=0.6pt] (0,1.1) arc[start angle=90, end angle=270, radius=1.1];
        \node[ line width=0.6pt, dashed, draw opacity=0.5] (a) at (-1.3,0){$\scriptstyle a$};
        \node[ line width=0.6pt, dashed, draw opacity=0.5] (a) at (0,-1.7){$\scriptstyle x$};
        \node[ line width=0.6pt, dashed, draw opacity=0.5] (a) at (-0.25,-0.7){$\scriptstyle y$};
        \node[ line width=0.6pt, dashed, draw opacity=0.5] (a) at (0.26,-1.15){$\scriptstyle \nu$};
        \node[ line width=0.6pt, dashed, draw opacity=0.5] (a) at (-0.25,0.7){$\scriptstyle z$};
        \node[ line width=0.6pt, dashed, draw opacity=0.5] (a) at (0,1.7){$\scriptstyle w$};
        \node[ line width=0.6pt, dashed, draw opacity=0.5] (a) at (0.2,1.15){$\scriptstyle \zeta$};
    \end{tikzpicture}
\end{aligned}
    \right)  \\
     = & \frac{\delta_{s,w}\delta_{t,x}\delta_{u,y}\delta_{v,z}}{d_s}\;\begin{aligned}
        \begin{tikzpicture}[scale=0.7]
        \path[black!60, fill=gray!15] (0,0) circle[radius=1.4];
             \draw[line width=.6pt,black] (0,-1.4)--(0,1.4);
             \draw[red, line width=0.6pt] (0,0.95) arc[start angle=90, end angle=270, radius=0.7];
             \draw[blue, line width=0.6pt] (0,-0.95) arc[start angle=-90, end angle=90, radius=0.7];
             \draw[line width=.6pt,black] (0,1.4) arc[start angle=90, end angle=-90, radius=1.4];
             \node[ line width=0.6pt, dashed, draw opacity=0.5] (a) at (1.2,0){$\scriptstyle \bar{s}$};
            \node[ line width=0.6pt, dashed, draw opacity=0.5] (a) at (-0.9,0.4){$\scriptstyle a$};
            \node[ line width=0.6pt, dashed, draw opacity=0.5] (a) at (0.9,-0.4){$\scriptstyle b$};
            \node[ line width=0.6pt, dashed, draw opacity=0.5] (a) at (-0.2,-0.7){$\scriptstyle t$};
            \node[ line width=0.6pt, dashed, draw opacity=0.5] (a) at (-0.2,-1.1){$\scriptstyle \mu$};
            \node[ line width=0.6pt, dashed, draw opacity=0.5] (a) at (0.2,-0.5){$\scriptstyle \nu$};
            \node[ line width=0.6pt, dashed, draw opacity=0.5] (a) at (-0.2,-0.1){$\scriptstyle u$};
            \node[ line width=0.6pt, dashed, draw opacity=0.5] (a) at (-0.2,0.7){$\scriptstyle v$};
            \node[ line width=0.6pt, dashed, draw opacity=0.5] (a) at (-0.2,0.3){$\scriptstyle \gamma$};
            \node[ line width=0.6pt, dashed, draw opacity=0.5] (a) at (0.2,1){$\scriptstyle \zeta$};
        \end{tikzpicture}
    \end{aligned}\;. 
\end{aligned}
\end{equation}
It was proved in \cite{jia2025weakhopftubealgebra} that this defines a skew-pairing. Furthermore, the right-hand side of Eq.~\eqref{eq:pairing_def} can be explicitly evaluated using topological local move to yield a complex number.
We will not distinguish $\hat{H}_3$ and $\mathbf{Tube}({\mathcal{H}_3}_{\mathcal{H}_3})^{\rm cop}$ anymore hereinafter. The Haar integral of $\hat{H}_3$ is given by
\begin{equation}\label{eq:HaarMTube}
    \Lambda = \frac{1}{\operatorname{rank} \eC} \sum_{a,x,y,\mu} 
    \sqrt{\frac{d_a}{d_x^3 d_y}}
        \begin{aligned}
    \begin{tikzpicture}
        \begin{scope}
            \fill[gray!20]
                (0,1.1) arc[start angle=90, end angle=-90, radius=1.1] -- 
                (0,-0.5) arc[start angle=-90, end angle=90, radius=0.5] -- cycle;
        \end{scope}
            \draw[line width=0.6pt,black,->] (0,0.5)--(0,0.7);
            \draw[line width=0.6pt,black,->] (0,0.7)--(0,1.0);
           \draw[line width=0.6pt,black] (0,0.5)--(0,1.1);
        \draw[line width=.6pt,black] (0,-0.5)--(0,-1.1);
        \draw[line width=0.6pt,black,->] (0,-1.1)--(0,-0.9);
        \draw[line width=0.6pt,black,->] (0,-1.1)--(0,-0.6);
        \draw[blue, line width=0.6pt] (0,-0.8) arc[start angle=-90, end angle=90, radius=0.8];
          \draw[blue, line width=0.6pt, ->] (0,-0.8) arc[start angle=-90, end angle=0, radius=0.8];
        \node[line width=0.6pt, dashed, draw opacity=0.5] at (0,1.3) {$y$};
        \node[line width=0.6pt, dashed, draw opacity=0.5] at (0,-1.3) {$y$};
        \node[line width=0.6pt, dashed, draw opacity=0.5] at (1,0) {$a$};
        \node[line width=0.6pt, dashed, draw opacity=0.5] at (-0.3,-0.7) {$\mu$};
        \node[line width=0.6pt, dashed, draw opacity=0.5] at (0,-0.3) {$x$};
        \node[line width=0.6pt, dashed, draw opacity=0.5] at (0,0.3) {$x$};
        \node[line width=0.6pt, dashed, draw opacity=0.5] at (-0.3,0.7) {$\mu$};
    \end{tikzpicture}
\end{aligned}.
\end{equation}
The Haar integral of $\hat{H}_3$ induces a inner product of $H_3$,
\begin{equation}
    \langle x,y\rangle :=\Lambda(x^*y).
\end{equation}
For this reason $\Lambda$ is also called of Haar measure for $H_3$.

\black

\vspace{0.5em}
\emph{Symmetry Topological Field Theory (SymTFT) with Haagerup $\cH_3$ symmetry.} --- SymTFT provides a general framework for studying non-invertible symmetry-protected topological (SPT) phases in both gapped and gapless cases (see, e.g.,~\cite{Kong2020algebraic, SchaferNameki2024ICTP, huang2023topologicalholo, bhardwaj2024lattice, freed2024topSymTFT, gaiotto2021orbifold, bhardwaj2023generalizedcharge, apruzzi2023symmetry, bhardwaj2024gappedphases, Zhang2024anomaly, Ji2020categoricalsym}).
The SymTFT can be applied to weak Hopf lattice gauge theory and establish the general algebraic framework for \((1+1)\)D non-invertible SPT phases \cite{jia2024generalized,jia2024weakhopf}. 
Based on this, we will construct a lattice model.

The SymTFT has a sandwich structure, as illustrated in Fig.~\ref{fig:ladder}~(a). The sandwich manifold is defined as \( \mathbb{M}^{1,1} \times [0,1] \), where \( \mathbb{M}^{1,1} \) represents the \((1+1)\)-dimensional manifold on which our system resides. This manifold has two boundaries: the basic idea is to impose non-invertible symmetry on one boundary and place the physical system on the other. After compactifying over the interval \( [0,1] \), we obtain a \((1+1)\)-dimensional system with non-invertible symmetry. To construct the SymTFT, the symmetry boundary must be chosen as a gapped topological boundary condition. There are many topological orders that do not have any gapped boundary, but non-chiral topological orders are guaranteed to have at least one gapped boundary. Thus, we put a non-chiral topological order on the 2d bulk.

For Haagerup symmetry \( \mathcal{H}_3 = \Rep(H_3) \), we place \( \mathcal{Z}(\cH_3) \) on the \((2+1)\)D bulk of SymTFT, which means the bulk is a weak Hopf quantum double model with input data given by \( H_3 \). The bulk topological excitation is characterized by \( \mathcal{Z}(\mathcal{H}_3) \simeq_{\otimes, \rm br} \Rep(D(H_3)) \), where the equivalence is both monoidal and braided, and \( D(H_3) \) denotes the quantum double of \( H_3 \) (see~\cite{Jia2023weak} for details on the weak Hopf quantum double model).
The symmetry boundary is chosen as \( \eB_{\rm sym} = \mathcal{H}_3 = \Rep(H_3) \), meaning that the symmetry boundary is a smooth boundary for the quantum double model. The boundary input data is \( H_3 \) (regarded as a comodule algebra over \( H_3 \)).
For the physical boundary, we can choose either a gapped or gapless boundary condition, and the resulting lattice model will possess Haagerup fusion category symmetry.

\vspace{0.5em}
\emph{Haagerup Fusion Symmetric Cluster Ladder Model.} --- 
To construct the lattice model with Haagerup fusion category symmetry, we first introduce the generalized Pauli operators. The regular action of a Haagerup weak Hopf algebra \( H_3 \) on itself can be viewed as a generalization of Pauli \( X \)-type operators. For the left action \( H_3 \curvearrowright H_3 \), we define:
\begin{equation}
\XR_g |h\rangle = |gh\rangle, \quad \XL_g |h\rangle = |hS^{-1}(g)\rangle, \quad \forall g,h \in H_3.
\end{equation}
The dual space \( \hat{H}_3 := \Hom(H_3, \mathbb{C}) \) is also a weak Hopf algebra. There are also canonical actions of the dual weak Hopf algebra \( \hat{H}_3 \) on the Hopf qudit \( H_3 \), defined using Sweedler’s notation as follows:
\begin{equation}
\varphi \rightharpoonup x := \sum_{(x)} x^{(1)} \langle \varphi, x^{(2)} \rangle, \quad x \leftharpoonup \varphi := \sum_{(x)} \langle \varphi, x^{(1)} \rangle x^{(2)},
\end{equation}
for all \( \varphi \in \hat{H}_3 \) and \( x \in H_3 \). For the left action \( \hat{H}_3 \curvearrowright H_3 \), we define:
\begin{equation}
\ZR_{\psi} |h\rangle = |\psi \rightharpoonup h\rangle = \sum_{(h)} \psi(h^{\ctwo}) |h^{\cone}\rangle,
\end{equation}
\begin{equation}
\ZL_{\psi} |h\rangle = |h \leftharpoonup \hat{S}(\psi)\rangle = \sum_{(h)} \psi(S(h^{\cone})) |h^{\ctwo}\rangle.
\end{equation}
These operators can be regarded as generalized Pauli \( Z \)-type operators.

Our construction of local operators is based on the comodule algebra over \( H_3 \). By definition, an algebra \( K \) is referred to as a left \( H_3 \)-comodule algebra~\cite{Bohmdoihopf} if there exists a map \( \beta: K \to H_3 \otimes K \) that satisfies the following properties:
\begin{equation}
    \beta(xy) = \beta(x) \beta(y),
\end{equation}
and
\begin{equation} \label{eq:weakComodulealge}
    (1_{H_3} \otimes x) \beta(1_K) = (\varepsilon_R \otimes \id_K) \circ \beta(x),
\end{equation}
for all \( x, y \in K \). 
Also see Ref. \cite{jia2024weakhopf} for several equivalent definitions.
Using Sweedler's notation, the coaction \( \beta(x) \) is expressed as \( \sum_{[x]} x^{[-1]} \otimes x^{[0]} \), with higher coactions represented as \( \beta_2(x) = \sum_{[x]} x^{[-2]} \otimes x^{[-1]} \otimes x^{[0]} \), and so forth. Each topological boundary condition corresponds to an indecomposable comodule algebra \( K \) over \( H_3 \). This comodule algebra also corresponds to an indecomposable module category \( \EuScript{M}_K \) over \( \cH_3 \), as well as a Lagrangian algebra \( \EuScript{A}_K \) in the bulk quantum double phase \( \mathcal{Z}(\cH_3) \).

The model is put on an ultra-thin sandwich lattice (regarded as a weak Hopf quantum double model \cite{Jia2023weak}):
\begin{equation}\label{eq:ClusterLadder}
\begin{aligned}
\begin{tikzpicture}
    \def\n{5}
    \def\s{1}
    
    \fill[yellow!10] (0, 0) rectangle (\n*\s+\s, \s); 

    \foreach \i in {0,...,\n} {
        \draw[-stealth, line width=1.0pt,blue, midway] (\i*\s, 0) -- (\i*\s+\s, 0);
        \draw[-stealth, line width=1.0pt,red, midway] (\i*\s, \s) -- (\i*\s+\s, \s);
        \draw[-stealth,line width=1.0pt, midway] (\i*\s, 0) -- (\i*\s, \s);
    }
    \draw[-stealth, midway,line width=1.0pt] (\n*\s+\s, 0) -- (\n*\s+\s, \s);
\end{tikzpicture}
\end{aligned}
\end{equation}
Periodic or open boundary conditions can be applied in the horizontal direction.
The bulk edges are drawn in black, with the edge Hilbert space chosen as the Haagerup weak Hopf qudit \( H_3 \). The edges of the symmetry boundary are drawn in blue, with the edge space given by the left comodule algebra \( K = H_3 \). The edges of the physical boundary are drawn in red, where the edge space can be chosen as an arbitrary right \( H_3 \)-comodule algebra \( J \).
This means that the total Hilbert space has well-defined tensor-product structure.

Observe that, in the described lattice structure, the bulk resides to the left of the symmetry boundary and to the right of the physical boundary when moving along the positive boundary direction. Therefore, the associated comodule algebras \( K \) and \( J \) should be selected as a left \( H_3 \)-comodule algebra and a right \( H_3 \)-comodule algebra, respectively. Modifying the configuration of edge orientations in the lattice leads to different models, but does not alter the underlying physics of the system.

Since our goal is to realize the Haagerup fusion category symmetry, we choose the symmetry boundary to be the smooth boundary, which corresponds to the left comodule algebra \( K = H_3 \).
In this case, the boundary vertex operator is given by:
\begin{equation}
    \begin{aligned}
     \begin{tikzpicture}
    \begin{scope}[rotate=90] 
        \fill[yellow!10] (0, 0) rectangle (1,2); 
        \draw[-stealth,blue,line width = 1.6pt] (0,1) -- (0,0); 
        \draw[-stealth,blue,line width = 1.6pt] (0,2) -- (0,1); 
        \draw[-stealth,black] (0,1) -- (1,1); 
        \node[line width=0.2pt, dashed, draw opacity=0.5] (a) at (-0.4,0.5) {\( y \)};
        \node[line width=0.2pt, dashed, draw opacity=0.5] (a) at (-0.4,1.5) {\( x \)};
        \node[line width=0.2pt, dashed, draw opacity=0.5] (a) at (0.5,0.8) {\( h \)};
        \node[line width=0.2pt, dashed, draw opacity=0.5] (a) at (0.2,0.5) {\textcircled{2}};
        \node[line width=0.2pt, dashed, draw opacity=0.5] (a) at (0.3,1.8) {\textcircled{1}};
        \node[line width=0.2pt, dashed, draw opacity=0.5] (a) at (0.8,1.2) {\textcircled{3}};
    \end{scope}
\end{tikzpicture}
\end{aligned} 
\quad 
\begin{aligned}
    & \Av^K_{v_s} = \XR_{\lambda^{\cone}} \otimes \XL_{\lambda_K^{\ctwo}} \otimes \XL_{\lambda^{\cthree}}, \\
    & \Av^K_{v_s}|x,y,h\rangle = |\lambda^{\cone} x, yS^{-1}(\lambda^{\ctwo}), hS^{-1}(\lambda^{\cthree})\rangle.
\end{aligned}
\label{eq:Avs}
\end{equation}
where \( \lambda \) is the Haar integral of \( H_3 \) (Eq.~\eqref{eq:HaarTube}).
The cocommutativity of \( \lambda \) ensures that the operator \( \Av^K_{v_s} \) is independent of the specific starting link \( s = (v_s, f) \) chosen at the vertex \( v_s \). Consequently, \( \Av^K_{v_s} \) depends only on the vertex \( v_s \) itself. The symmetry boundary Hamiltonian is defined as
\begin{equation}
  \mathbb{H}_{\rm sym} = -\sum_{v_s} \Av^K_{v_s},
\end{equation}
where the sum runs over all vertices in the symmetry boundary.

The face operator is constructed from generalized Pauli \( Z \) operators as follows:
\begin{equation}
\begin{aligned}
    \begin{tikzpicture}
    \begin{scope}[rotate=-90] 
        \fill[yellow!10] (0, 0) rectangle (1,1); 
        \draw[-stealth,red,line width = 1.6pt] (0,0) -- (0,1);
        \draw[cyan,dotted, line width = 1pt] (0,0) -- (0.5,0.5); 
        \draw[-stealth,blue,line width = 1.6pt] (1,0) -- (1,1);
        \draw[-stealth,black] (1,1) -- (0,1); 
        \draw[-stealth,black] (1,0) -- (0,0); 
        \node[line width=0.2pt, dashed, draw opacity=0.5] (a) at (-0.3,0.5) {\( x \)};
        \node[line width=0.2pt, dashed, draw opacity=0.5] (a) at (1.4,0.5) {\( y \)};
        \node[line width=0.2pt, dashed, draw opacity=0.5] (a) at (0.5,1.25) {\( h \)};
        \node[line width=0.2pt, dashed, draw opacity=0.5] (a) at (0.5,-0.25) {\( g \)};
        \node[line width=0.2pt, dashed, draw opacity=0.5] (a) at (0.25,0.5) {\textcircled{4}};
        \node[line width=0.2pt, dashed, draw opacity=0.5] (a) at (0.8,0.5) {\textcircled{2}};
        \node[line width=0.2pt, dashed, draw opacity=0.5] (a) at (0.5,0.8) {\textcircled{3}};
        \node[line width=0.2pt, dashed, draw opacity=0.5] (a) at (0.5,0.2) {\textcircled{1}};
    \end{scope}
\end{tikzpicture}   
\end{aligned}
\quad
\begin{aligned}
   & \Bf_f^{\psi} = \ZR_{\psi^{\cone}} \otimes \ZL^K_{\psi^{\ctwo}} \otimes \ZL_{\psi^{\cthree}} \otimes \ZR_{\psi^{\cfour}}^J.
\end{aligned}
\end{equation}
Written explicitly, we have
\[
\Bf_f^{\psi} |g,y,h,x\rangle = \sum \psi\left(g^{\ctwo} S(y^{[-1]}) S(h^{(1)}) x^{[1]}\right) |g^{\cone},y^{[0]},h^{\ctwo},x^{[0]}\rangle.
\]
By setting \( \psi \) to be the Haar measure \( \Lambda \in \hat{H}_3 \), we obtain the face stabilizer operator \( \Bf_f = \Bf_f^{\Lambda} \), and the corresponding bulk Hamiltonian is
\begin{equation}
    \Hbb_{\rm bk} = -\sum_f \Bf_f,
\end{equation}
where the sum runs over all faces of the lattice.

The physical boundary can be chosen to be either gapped or gapless. In the gapped case, we assign a right \( H_3 \)-comodule algebra to each edge of the physical boundary. The vertex operator can be constructed in a similar way as in Eq.~\eqref{eq:Avs} (see Ref.~\cite{jia2024weakhopf}, where we need to use the symmetric separability idempotent of \( J \)). The physical boundary Hamiltonian is given by
\begin{equation}
\Hbb_{\rm phys} = -\sum_{v_p} \Av_{v_p}^J,
\end{equation}
where the sum runs over all vertices \( v_p \) of the physical boundary.


The cluster ladder model exhibiting Haagerup fusion category symmetry is given by
\begin{equation}\label{eq:SymHam}
    \mathbb{H}_{H_3}^{K,J} = \mathbb{H}_{\rm bk} + \mathbb{H}_{\rm sym} + \mathbb{H}_{\rm phys}.
\end{equation}
It can be shown that all local terms in the Hamiltonian mutually commute. Fixing the symmetry boundary to be the smooth boundary, i.e., \( K = H_3 \), different choices of comodule algebras \( J \) over the Haagerup weak Hopf algebra \( H_3 \) lead to distinct \( \mathcal{H}_3 \)-symmetric phases.
In addition to imposing a topological boundary condition on the physical boundary, one can also consider non-topological or even gapless boundary conditions, which can give rise to more intricate quantum phases with \( \mathcal{H}_3 \) symmetry.

In the algebraic theory of weak Hopf SymTFT, the topological excitations on the symmetry boundary are characterized by $H_3$-covariant $K|K$-bimodules, denoted by ${}_K^{H_3}\mathsf{Mod}_K$. Similarly, for the physical boundary characterized by the comodule algebra $J$, the corresponding topological excitations are described by $H_3$-covariant $J|J$-bimodules, ${}_J\mathsf{Mod}_J^{H_3}$.
To relate this to the string-net picture, we consider the bulk fusion category $\mathsf{Rep}(H_3) \simeq \mathcal{H}_3$. Both the symmetry and physical boundaries are described by module categories over $\mathcal{H}_3$. Specifically, the category of left $K$-modules, $_K\mathsf{Mod}$, forms a left module category over $\mathcal{H}_3$, while the category of right $J$-modules, $\mathsf{Mod}_J$, forms a right module category over $\mathcal{H}_3$.
The distinction between left and right module categories is determined by the orientation of the lattice edges. In this work, we fix all boundary edges to be oriented rightwards, and this choice sets our convention.
See Table~\ref{tab:SNdic} for a summary.

\begin{table}[h]
        \centering
\caption{Dictionary between the data of the $D(H_3)$ quantum double model and the $\mathcal{H}_3$ string-net model.\label{tab:SNdic}}  
\resizebox{1.0\linewidth}{!}{%
\begin{tabular}{l|l|l}
        \hline
                \hline
           & QD model (quantum group)  & String-net (fusion category) \\
             \hline
       bulk    & weak Hopf algebra $H_3$ &  fusion category $\mathsf{Rep}(H_3)\simeq \mathcal{H}_3$\\
       \hline
    boundary    &  comodule algebra $K$ & module category  $\mathcal{M}\simeq \mathsf{Mod}_K$\\
       \hline
    bd. ext.    & $H_3$-covariant $K|K$-bimodule & module functor  $\mathsf{Fun}_{\mathcal{H}_3}(\mathcal{M},\mathcal{M})$\\
       \hline
            \hline
\end{tabular}
}
\end{table}

In fact, the model possesses a larger symmetry than $\mathcal{H}_3$~\cite{jia2024weakhopf}. On a closed manifold, the symmetry arising from the symmetry boundary is given by the cocommutative subalgebra $\Cocom(\hat{H}_3)$ of the dual weak Hopf algebra $\hat{H}_3$; on an open manifold, the symmetry can be as large as the full dual algebra $\hat{H}_3$ itself.
Since \( \mathcal{H}_3 \simeq \Rep(H_3) \), the Grothendieck group \( \operatorname{Gr}(\mathcal{H}_3) \) is generated by the irreducible characters \( \chi_{\Gamma} \) of \( H_3 \), with multiplication rule
\begin{equation}
\chi_{\Gamma} \cdot \chi_{\Phi} = \chi_{\Gamma \otimes \Phi} = \sum_{\Psi} N_{\Gamma, \Phi}^{\Psi} \chi_{\Psi},
\end{equation}
where \( N_{\Gamma, \Phi}^{\Psi} \in \mathbb{Z}_{\geq 0} \) are the fusion multiplicities. The resulting character algebra \( R(\mathcal{H}_3) := \operatorname{Gr}(\mathcal{H}_3) \otimes_{\mathbb{Z}} \mathbb{C} \) forms a subalgebra of the dual weak Hopf algebra \( \hat{H}_3 \), all elements in it are cocommutative since all character functions are cocommutative.
This implies that the Haagerup fusion category symmetry is a sub-symmetry of the dual weak Hopf symmetry \( \hat{H}_3 \) on open manifolds, and of its cocommutative subalgebra \( \Cocom(\hat{H}_3) \) on closed manifolds.

\vspace{1em}
\emph{Cluster state model as quantum double model with smooth and rough boundary.} ---
A crucial special case is the cluster state model, which can be viewed as an ultra-thin quantum double model with one rough and one smooth boundary~\cite{jia2024generalized,jia2024weakhopf}. We emphasize that, in this context, the notion of a rough boundary refers to the removal of degrees of freedom at the boundary, a convention that has been commonly used in lattice models since the original work of Ref.~\cite{bravyi1998quantum}. 
This interpretation differs from another common usage in the literature, where a rough boundary corresponds to the module category $\mathsf{Vect}$ over the input fusion category in the string-net model framework, or equivalently, to the comodule algebra $\mathbb{C}$ over the input weak Hopf algebra in the quantum double model framework. Since $\mathcal{H}_3$ exhibits anomalous symmetry, $\mathsf{Vect}$ is not a module category over $\mathcal{H}_3$, and $\mathbb{C}$ is not a comodule algebra over $H_3$. However, for anomaly-free symmetries, these two notions coincide.
For cluster state, the lattice becomes (we assume the horizontal periodic boundary condition)
\begin{equation*}
\begin{aligned}
\begin{tikzpicture}
    \def\n{5}
    \def\s{1}   
    \fill[yellow!10] (0, 0) rectangle (\n*\s+\s, \s);

    \foreach \i in {0,...,\n} {
        \draw[-stealth, line width=1pt, blue, midway] (\i*\s, 0) -- (\i*\s+\s, 0);
        \draw[-stealth, line width=1pt, midway] (\i*\s, 0) -- (\i*\s, \s);
    }
    \draw[-stealth, line width=1pt, midway] (\n*\s+\s, 0) -- (\n*\s+\s, \s);
    \draw[-stealth, white, line width=2pt, midway] (6, 0) -- (6, 1.02);
\end{tikzpicture}
\end{aligned}
\end{equation*}
The Hamiltonian is thus
\begin{equation}\label{eq:clusterSH}
    \mathbb{H}_{\rm cluster}=-\sum_{v_s} \Av_{v_s} -\sum_{f} \Bf_f.
\end{equation}
This is analogous to the CSS-type cluster state model, which features only two types of local stabilizers: vertex operators and face operators~\cite{brell2015generalized,fechisin2023noninvertible,jia2024generalized} (see~\cite{jia2024generalized} for a detailed discussion of this perspective).

To find the symmetry of the Haagerup cluster state model, we deform the ladder (with periodic boundary conditions) into a cone. We introduce a vertex \( v_{\rm rough} \) on the rough boundary (the apex of the cone) and a face \( f_{\rm smooth} \) on the smooth boundary (the base of the cone):
\begin{equation*}
    \begin{tikzpicture}
    \fill[red!60] (0,0) ellipse (1.5 and 0.5); 
   
    \draw[blue!60, very thick, -latex] (1.5,0) arc[start angle=0, end angle=60, x radius=1.5, y radius=0.5];

    \draw[blue!60, very thick, -latex] (1.5,0) arc[start angle=0, end angle=130, x radius=1.5, y radius=0.5];
    \draw[blue!60, very thick, -latex] (1.5,0) arc[start angle=0, end angle=220, x radius=1.5, y radius=0.5];
    \draw[blue!60, very thick, -latex] (1.5,0) arc[start angle=0, end angle=320, x radius=1.5, y radius=0.5];
    \draw[blue!60, very thick] (0,0) ellipse (1.5 and 0.5); 

    \coordinate (Apex) at (0,1.5); 

    \fill[yellow!10, opacity=0.7] (Apex) -- (-1.5,0) arc[start angle=180, end angle=360, x radius=1.5, y radius=0.5] -- cycle;

    \draw[thick] (1.52,0) -- node[midway] {} (Apex);  
    \draw[thick] (-1.52,0) -- node[midway] {} (Apex); 
    \draw[thick] (0.3,0.45) -- node[midway] {} (Apex);   
    \draw[thick] (-0.3,-0.53) -- node[midway] {} (Apex);  

    \coordinate (RightMid) at (1.52,0); 
    \draw[thick,  -latex] (RightMid) -- node[midway, left] {} ($ (Apex)!0.5!(1.52,0) $);  

    \coordinate (LeftMid) at (-1.52,0); 
    \draw[thick,  -latex] (LeftMid) -- node[midway, left] {} ($ (Apex)!0.5!(-1.52,0) $);  

    \coordinate (TopMid) at (0.3,0.45); 
    \draw[thick,  -latex] (TopMid) -- node[midway, right] {} ($ (Apex)!0.5!(0.3,0.45) $);  

    \coordinate (BottomMid) at (-0.3,-0.53); 
    \draw[thick,  -latex] (BottomMid) -- node[midway, right] {} ($ (Apex)!0.5!(-0.3,-0.53) $);  

    \draw[cyan, dotted] (0,0) -- (1.52,0);  
    \draw[cyan, dotted] (0,0) -- (-1.52,0); 
    \draw[cyan, dotted] (0,0) -- (0.3,0.53);  
    \draw[cyan, dotted] (0,0) -- (-0.3,-0.53);  
    \end{tikzpicture}
\end{equation*}
In this way, we can regard the model as a quantum double model on a sphere with one face operator and one vertex operator removed.

The symmetry operator at the smooth boundary is the face operator \( W_{\varphi} = \Bf^{\varphi}_{f_{\rm smooth}} \). 
By definition, this operator takes the form
\begin{equation} 
W_{\varphi} = \sum_{(\varphi)} \ZR_{\varphi^{(1)}} \otimes \ZR_{\varphi^{(2)}} \otimes \cdots \otimes \ZR_{\varphi^{(n)}},
\end{equation}
where $\Delta^{(n-1)}(\varphi) = \sum_{(\varphi)} \varphi^{(1)} \otimes \varphi^{(2)} \otimes \cdots \otimes \varphi^{(n)}$ is the $(n{-}1)$-fold coproduct in $\hat{H}_3$. More precisely, we have
\begin{equation} \label{eq:symOPZ}
\begin{aligned}
     \begin{tikzpicture}
   \fill[yellow!10] (0, 0) rectangle (5,1); 
    \draw[-stealth, line width=1.0pt, blue, midway] (0,0) --(1,0);
    \draw[-stealth, line width=1.0pt, blue, midway] (1,0) --(2,0);
    \draw[-stealth, line width=1.0pt, blue, midway] (2,0) --(2.5,0);
    \draw[-stealth, line width=1.0pt, blue, midway] (3.5,0) --(4,0);
    \draw[-stealth, line width=1.0pt, blue, midway] (4,0) --(5,0);
      \draw[-stealth,line width=1.0pt, midway] (0, 0) -- (0, 1);
     \draw[-stealth,line width=1.0pt, midway] (1, 0) -- (1, 1); 
    \draw[-stealth,line width=1.0pt, midway] (2, 0) -- (2, 1); 
    \draw[-stealth,line width=1.0pt, midway] (4, 0) -- (4, 1); 
     \node[line width=0.2pt, dashed, blue, draw opacity=0.5] (a) at (0.5,-0.4) {\(\ZR_{\varphi^{\cone}} \)};
    \node[line width=0.2pt, dashed, blue, draw opacity=0.5] (a) at (1.5,-0.4) {\(\ZR_{\varphi^{\ctwo}} \)};
     \node[line width=0.2pt, dashed, blue, draw opacity=0.5] (a) at (3,-0.4) {\(\cdots \)};
    \node[line width=0.2pt, dashed, blue, draw opacity=0.5] (a) at (4.5,-0.4) {\(\ZR_{\varphi^{(n)}} \)};
    \end{tikzpicture}
\end{aligned}
\end{equation}
We have the relation $W_{\psi} W_{\phi} = W_{\psi \cdot \phi}$ for $\psi, \phi \in \hat{H}_3$. Note that for the unit $\hat{1}$ in $\hat{H}_3$, $W_{\hat{1}}$ is not equivalent to the identity operator, since $\Delta(\hat{1}) \neq \hat{1} \otimes \hat{1}$ in this case. However, we still have $W_{1}W_{\psi} = W_{\psi} = W_{\psi}W_1$. When $\varphi \in \Cocom(\hat{H}_3)$, the operator $\Bf^{\varphi}_{f_{\rm smooth}}$ commutes with the Hamiltonian~\cite{chang2014kitaev, Jia2023weak}. This implies that the model possesses a $\Cocom(\hat{H}_3)$ symmetry.
For \( \cH_3 = \Rep(H_3) \), its character algebra is contained in \( \Cocom(\hat{H}_3) \), and thus the model also has Haagerup fusion category symmetry (more precisely Haagerup fusion algebra symmetry).

The rough boundary also gives symmetry operators \( W_h = \Av_{v_{\rm rough}}^h \), with the relation \( W_h W_g = W_{hg} \). 
 More precisely, \( W_h = \sum_{(h)} \XR_{h^{\cone}} \otimes \cdots \otimes \XR_{h^{(n)}} \):
\begin{equation}
\begin{aligned}
     \begin{tikzpicture}
   \fill[yellow!10] (0, 0) rectangle (5,1); 
    \draw[-stealth, line width=1.0pt, blue, midway] (0,0) --(1,0);
    \draw[-stealth, line width=1.0pt, blue, midway] (1,0) --(2,0);
    \draw[-stealth, line width=1.0pt, blue, midway] (2,0) --(2.5,0);
    \draw[-stealth, line width=1.0pt, blue, midway] (3.5,0) --(4,0);
    \draw[-stealth, line width=1.0pt, blue, midway] (4,0) --(5,0);
      \draw[-stealth,line width=1.0pt, midway] (0, 0) -- (0, 1);
     \draw[-stealth,line width=1.0pt, midway] (1, 0) -- (1, 1); 
    \draw[-stealth,line width=1.0pt, midway] (2, 0) -- (2, 1); 
    \draw[-stealth,line width=1.0pt, midway] (4, 0) -- (4, 1); 
     \node[line width=0.2pt, dashed, draw opacity=0.5] (a) at (0.5,0.5) {\(\XR_{h^{\cone}} \)};
    \node[line width=0.2pt, dashed, draw opacity=0.5] (a) at (1.5,0.5) {\(\XR_{h^{\ctwo}} \)};
    \node[line width=0.2pt, dashed, draw opacity=0.5] (a) at (2.5,0.5) {\(\XR_{h^{\cthree}} \)};
     \node[line width=0.2pt, dashed, draw opacity=0.5] (a) at (3.4,0.5) {\(\cdots \)};
    \node[line width=0.2pt, dashed, draw opacity=0.5] (a) at (4.5,0.5) {\(\XR_{h^{(n)}} \)};
    \end{tikzpicture}
\end{aligned}
\end{equation}
When \( h \in \Cocom(H_3) \), \( W_h \) commutes with the Hamiltonian. This implies that the model has a \( \Cocom(H_3) \) symmetry.
For open manifolds, there is no need for the symmetry elements to be cocommutative, and thus the symmetry becomes \( \hat{H}_3 \times H_3 \).
The result is summarized in Theorem~\ref{thm:cluster}.

For the general cluster ladder model, due to the existence of boundary data \( J \), the symmetry \( H_3 \) from the rough boundary is broken. Using a similar discussion as above, we find that the dual symmetry \( \hat{H}_3 \) still exists. Thus, for periodic boundary conditions, the model has \( \Cocom(\hat{H}_3) \) symmetry, while for open boundary conditions, the model has \( \hat{H}_3 \) symmetry. In both cases, the symmetry contains \( \cH_3 \) as a sub-symmetry.

\vspace{0.5em}
\emph{Matrix-Product State Representation of Haagerup Cluster State.} --- 
The cluster state model can be solved using the weak Hopf matrix product state formalism developed in Refs.~\cite{jia2024generalized,jia2024weakhopf,Jia2023weak,jia2023boundary}. In this section, we will focus on the ground state for the cluster state model. For the more general cluster ladder model, the solution can be obtained following the discussion in Ref.~\cite{jia2024weakhopf}.

The main tool we employ is the comultiplication and canonical pairing between $H_3$ and $\hat{H}_3$, which satisfy the properties given in Eqs.~\eqref{eq:MultiplicationBarA}--\eqref{eq:AntipodeBarA}. The gluing of local tensors in $H_3$ and $\hat{H}_3$ is defined via this pairing. For instance, given $\psi \in \hat{H}_3$ and $h \in H_3$, we can express them in dual orthonormal bases $\langle \hat{e}_j,e_i\rangle=\hat{e}_j(e_i)=\delta_{ij}$:
\[ 
h = \sum_i h_i e_i, \quad \psi = \sum_{i} \psi_i \hat{e}_i.
\]
The pairing between them is then given by
\begin{equation}
    \langle \psi, h \rangle = \sum_{i,j} \langle \psi_j \hat{e}_j, h_i e_i \rangle = \sum_{i} \psi_i h_i.
\end{equation}
Thus, this reduces to the standard tensor contraction. The comultiplication $\Delta_n(h)$ (the $n$-fold comultiplication) generates $n+1$ legs for a tensor. Multiplication induces the gluing of tensors in $H_3$ via the structure constants: $e_i \cdot e_j = \sum_k A_{ij}^k e_k$, so that
$
hg = \sum_{i,j,k} h_i g_j A_{ij}^k e_k.
$
For further details, see Ref.~\cite{jia2024generalized}.

We derive the ground state for the weak Hopf quantum double model in Ref.~\cite{Jia2023weak}. To apply this to the cluster state model, we introduce the following three local tensors.
For edges on the smooth boundary, the local tensor is chosen as the comultiplication of the Haar integral \( \lambda \in H_3 \):
\begin{equation}
     \Delta(\lambda) = 
     \begin{aligned}
\begin{tikzpicture}
    \node[draw, fill=blue!20, minimum width=0.6cm, minimum height=0.6cm] (center) at (0,0) {$\scriptstyle\lambda$};
    \draw[line width=1.0pt] (center.north) -- ++(0,.5) node[above] {$\scriptstyle\lambda^{\ctwo}$};
    \draw[line width=1.0pt] (center.south) -- ++(0,-.5) node[below] {$\scriptstyle \lambda^{\cone}$};
\end{tikzpicture}
\end{aligned}.
\end{equation}
The local tensor for a bulk edge corresponds to \( \Delta_2(\lambda) \):
\begin{equation}
     \Delta(\lambda)= \begin{aligned}			
\begin{tikzpicture}
    \draw[line width=1.0pt] (center.north) (0,0) -- ++(0.8,0) node[right] {$\scriptstyle\lambda^{\ctwo}$};
    \draw[line width=1.0pt] (center.north) (0,0) -- ++(-0.8,0) node[left] {$\scriptstyle\lambda^{\cthree}$};
    \draw[line width=1.0pt ] (center.south) -- ++(0,-.5) node[below] {$\scriptstyle \lambda^{\cone}$};
    \node[draw, fill=gray!20, minimum width=0.6cm, minimum height=0.6cm] (center) at (0,0) {$\scriptstyle\lambda$};
\end{tikzpicture}
\end{aligned}.
\end{equation}
The local tensor for a face corresponds to the Haar measure \( \Lambda \in \hat{H}_3 \), which is introduced to glue the edge tensors together:
\begin{equation}
    ( \id \otimes \id\otimes \hat{S}) \circ \hat{\Delta}(\Lambda)=
    \begin{aligned}			
\begin{tikzpicture}
\draw[line width=1.0pt] (center.north) (0,0) -- ++(0.8,0) node[right] {$\scriptstyle\Lambda^{\ctwo}$};
\draw[line width=1.0pt] (center.north) (0,0) -- ++(-0.8,0) node[left] {$\scriptstyle \hat{S}(\Lambda^{\cthree})$};
\draw[line width=1.0pt ] (center.south) -- ++(0,-.5) node[below] {$\scriptstyle \Lambda^{\cone}$};
\node[draw, fill=yellow!20, minimum size=0.6cm, shape=circle] (center) at (0,0) {$\scriptstyle\Lambda$};
 \filldraw[black] (-0.5,0) circle (2pt);  
\end{tikzpicture}
\end{aligned}. 
\end{equation}
The Haagerup cluster state is defined as the following tensor network state:
\begin{equation}
    \begin{aligned}			
\begin{tikzpicture}
    \draw[line width=1.0pt] (-0.5,0) -- ++(7.3,0) ;
    \draw[line width=1.0pt, red] (0,0) -- ++(0,-.8);
        \draw[line width=1.0pt, red] (2,0) -- ++(0,-.8);
    \draw[line width=1.0pt, red] (4,0) -- ++(0,-.8);
    \draw[line width=1.0pt, red] (6,0) -- ++(0,-.8);
     \filldraw[black] (0.5,0) circle (2pt);  
      \filldraw[black] (2.5,0) circle (2pt);
       \filldraw[black] (4.5,0) circle (2pt);
        \filldraw[black] (6.5,0) circle (2pt);
    \draw[line width=1.0pt] (1,-1) -- ++(0,1) ;
    \draw[line width=1.0pt, red] (1,-1) -- ++(0,-.8);
    \node[draw, fill=blue!20, minimum width=0.6cm, minimum height=0.6cm] (center) at (1,-1) {$\scriptstyle\lambda$};
     \draw[line width=1.0pt] (3,-1) -- ++(0,1) ;
    \draw[line width=1.0pt, red] (3,-1) -- ++(0,-.8);
    \node[draw, fill=blue!20, minimum width=0.6cm, minimum height=0.6cm] (center) at (3,-1) {$\scriptstyle\lambda$};
         \draw[line width=1.0pt] (5,-1) -- ++(0,1) ;
    \draw[line width=1.0pt, red] (5,-1) -- ++(0,-.8);
    \node[draw, fill=blue!20, minimum width=0.6cm, minimum height=0.6cm] (center) at (5,-1) {$\scriptstyle\lambda$};
    \node[draw, fill=gray!20, minimum width=0.6cm, minimum height=0.6cm] (center) at (0,0) {$\scriptstyle\lambda$};
    \node[draw, fill=yellow!20, minimum size=0.6cm, shape=circle] (center) at (1,0) {$\scriptstyle\Lambda$};
     \node[draw, fill=gray!20, minimum width=0.6cm, minimum height=0.6cm] (center) at (2,0) {$\scriptstyle\lambda$};
      \node[draw, fill=yellow!20, minimum size=0.6cm, shape=circle] (center) at (3,0) {$\scriptstyle\Lambda$};
        \node[draw, fill=gray!20, minimum width=0.6cm, minimum height=0.6cm] (center) at (4,0) {$\scriptstyle\lambda$};
      \node[draw, fill=yellow!20, minimum size=0.6cm, shape=circle] (center) at (5,0) {$\scriptstyle\Lambda$};
        \node[draw, fill=gray!20, minimum width=0.6cm, minimum height=0.6cm] (center) at (6,0) {$\scriptstyle\lambda$};
\end{tikzpicture}
\end{aligned} 
\end{equation}
The red free leg represents the physical degrees of freedom (which take values in $H_3$). This tensor network forms the ground state of the Haagerup cluster state model. To verify that this is indeed the ground state of Eq.~\eqref{eq:clusterSH}, we need to show that both $\Av_{v_s}$ and $\Bf_f$ stabilize the state. This is rigorously proved in \cite[Theorem 4]{Jia2023weak}, based on the properties $\lambda^2 = \lambda$ and $\Lambda^2 = \Lambda$, as well as the fact that both $\lambda$ and $\Lambda$ are cocommutative.

Recall that Haagerup symmetry can be expressed  as generalized Pauli-Z operators $\ZR$ via using the character $\chi_{\Gamma}$ of irreps $\Gamma$ of $H_3$ as in Eq.~\eqref{eq:symOPZ}.
It is convenient to introduce 
\begin{equation}
    \ZR_{\Gamma}|h\rangle =\sum_{(h)}|h^{\cone}\rangle \otimes \Gamma_{\alpha\beta}(h^{\ctwo}).
\end{equation}
If we take the trace over representation space (denote $\Tr'$) we obtain 
\begin{equation}
    \ZR_{\chi_{\Gamma}} = \Tr'  \ZR_{\Gamma}.
\end{equation}
This matches with the definition of $\ZR_{\chi_{\Gamma}}=(\id \otimes \chi_{\Gamma})\circ\Delta$ given in Eq.~\eqref{eq:symOPZ}. 
We can represent the comultiplication $\Delta$ using structure constants as $\Delta(e_i) = \sum_{j,k} C_{i}^{jk} e_j \otimes e_k$, which can be viewed as a tensor (read downward).
\begin{equation}
\Delta =
\begin{aligned}
    \begin{tikzpicture}
 \path[draw, fill=red!20]
      (-0.3,0) -- (0.3,0) -- (0,0.52) -- cycle; 
    \node at (0,0.17) {$\scriptstyle\Delta$};
    \draw[line width=1.0pt] (0,.52) -- ++(0,.3) node[above] {$\scriptstyle h$};
    \draw[line width=1.0pt] (-0.2,0) -- ++(0,-.3) node[below] {$\scriptstyle h^{\cone}$};
    \draw[line width=1.0pt] (0.2,0) -- ++(0,-.3) node[below] {$\scriptstyle h^{\ctwo}$};
\end{tikzpicture}
\end{aligned}.
\end{equation}
We can add more bottom legs to represent the $n$-fold comultiplication $\Delta_n$. The operator $\ZR_{\Gamma}$ can be represented as
\begin{equation}
     \ZR_{\Gamma} = 
     \begin{aligned}
\begin{tikzpicture}
    \draw[dotted,line width=1.0pt] (-.8,0) -- ++(1.6,0);
    \node[draw, fill=green!20, minimum width=0.6cm, minimum height=0.6cm] (center) at (0,0) {$\scriptstyle\ZR_{\Gamma}$};
    \draw[line width=1.0pt] (center.north) -- ++(0,.5) node[above] {$\scriptstyle $};
    \draw[line width=1.0pt] (center.south) -- ++(0,-.5) node[below] {$\scriptstyle $};
\end{tikzpicture}
\end{aligned}= \begin{aligned}
    \begin{tikzpicture}
 \path[draw, fill=red!20]
      (-0.3,0) -- (0.3,0) -- (0,0.52) -- cycle; 
    \node at (0,0.17) {$\scriptstyle\Delta$};
    \draw[line width=1.0pt] (0,.52) -- ++(0,.3) node[above] {$\scriptstyle h$};
    \draw[line width=1.0pt] (-0.2,0) -- ++(0,-.9);
    \draw[line width=1.0pt] (0.2,0) -- ++(0,-.9) ;
        \draw[dotted,line width=1.0pt] (-.5,-0.5) -- (0.9,-0.5);
    \path[draw, fill=cyan!20]
      (0,-0.3) rectangle ++(0.4,-0.4) ; 
          \node at (0.17,-0.5) {$\scriptstyle\Gamma$};
\end{tikzpicture}
\end{aligned}
\end{equation}
where we use dotted legs to represent the degrees of freedom in the representation space.
We can represent $\ZR_{\chi_{\Gamma}}$ as the result of contracting the dotted leg in the $\ZR_{\Gamma}$ tensor. Notice that
$
\left\langle \sum_{(\chi)} \chi^{(1)} \otimes \cdots \otimes \chi^{(n)},\, g_1 \otimes \cdots \otimes g_n \right\rangle = \chi_{\Gamma}(g_1 \cdots g_n),
$
which can be expressed as
\begin{equation}
  \ZR_{\chi_{\Gamma}^{(1)}} \otimes \cdots \otimes \ZR_{\chi_{\Gamma}^{(n)}} = \mathrm{Tr}'\left[\ZR_{\Gamma}(1) \star \cdots \star \ZR_{\Gamma}(n)\right],
\end{equation}
where we use $\star$ to denote multiplication in the representation space.
In this way, the $\mathcal{H}_3\simeq \Rep(H_3)$ symmetry operator can be expressed as 
\begin{equation}
    W_{\Gamma} =\mathrm{Tr}'\left[\ZR_{\Gamma}(1) \star \cdots \star \ZR_{\Gamma}(n)\right].
\end{equation}
When acting on cluster state, we have 
\begin{equation}
    \begin{aligned}			
\begin{tikzpicture}
    \draw[line width=1.0pt] (-0.5,0) -- ++(7.3,0) ;
    \draw[line width=1.0pt, red] (0,0) -- ++(0,-.8);
        \draw[line width=1.0pt, red] (2,0) -- ++(0,-.8);
    \draw[line width=1.0pt, red] (4,0) -- ++(0,-.8);
    \draw[line width=1.0pt, red] (6,0) -- ++(0,-.8);
     \filldraw[black] (0.5,0) circle (2pt);  
      \filldraw[black] (2.5,0) circle (2pt);
       \filldraw[black] (4.5,0) circle (2pt);
        \filldraw[black] (6.5,0) circle (2pt);
    \draw[line width=1.0pt] (1,-1) -- ++(0,1) ;
    \draw[line width=1.0pt, red] (1,-1) -- ++(0,-1.8);
    \node[draw, fill=blue!20, minimum width=0.6cm, minimum height=0.6cm] (center) at (1,-1) {$\scriptstyle\lambda$};
     \draw[line width=1.0pt] (3,-1) -- ++(0,1) ;
    \draw[line width=1.0pt, red] (3,-1) -- ++(0,-1.8);
    \node[draw, fill=blue!20, minimum width=0.6cm, minimum height=0.6cm] (center) at (3,-1) {$\scriptstyle\lambda$};
         \draw[line width=1.0pt] (5,-1) -- ++(0,1) ;
    \draw[line width=1.0pt, red] (5,-1) -- ++(0,-1.8);
    \node[draw, fill=blue!20, minimum width=0.6cm, minimum height=0.6cm] (center) at (5,-1) {$\scriptstyle\lambda$};
    \node[draw, fill=gray!20, minimum width=0.6cm, minimum height=0.6cm] (center) at (0,0) {$\scriptstyle\lambda$};
    \node[draw, fill=yellow!20, minimum size=0.6cm, shape=circle] (center) at (1,0) {$\scriptstyle\Lambda$};
     \node[draw, fill=gray!20, minimum width=0.6cm, minimum height=0.6cm] (center) at (2,0) {$\scriptstyle\lambda$};
      \node[draw, fill=yellow!20, minimum size=0.6cm, shape=circle] (center) at (3,0) {$\scriptstyle\Lambda$};
        \node[draw, fill=gray!20, minimum width=0.6cm, minimum height=0.6cm] (center) at (4,0) {$\scriptstyle\lambda$};
      \node[draw, fill=yellow!20, minimum size=0.6cm, shape=circle] (center) at (5,0) {$\scriptstyle\Lambda$};
        \node[draw, fill=gray!20, minimum width=0.6cm, minimum height=0.6cm] (center) at (6,0) {$\scriptstyle\lambda$};
       \draw[dotted,line width=1.0pt] (-.2,-2) -- ++(7,0);
    \node[draw, fill=green!20, minimum width=0.6cm, minimum height=0.6cm] (center) at (1,-2) {$\scriptstyle\ZR_{\Gamma}$};
     \node[draw, fill=green!20, minimum width=0.6cm, minimum height=0.6cm] (center) at (3,-2) {$\scriptstyle\ZR_{\Gamma}$};
     \node[draw, fill=green!20, minimum width=0.6cm, minimum height=0.6cm] (center) at (5,-2) {$\scriptstyle\ZR_{\Gamma}$};
\end{tikzpicture}
\end{aligned} 
\end{equation}
where the contraction is taken in the representation space.
We refer the reader to Refs.~\cite{jia2024generalized,jia2024weakhopf} for more details on weak Hopf tensor networks.

\vspace{0.5em}
\emph{Conclusion and Discussion.}---
In this work, we have proposed a cluster state model that exhibits Haagerup fusion category symmetry, and discussed its ground state and symmetry properties in detail. Despite this progress, several directions remain open for further exploration:

(i) \emph{Phases of the Haagerup cluster state model.} The phases of the model warrant deeper investigation. Since all structure constants can be calculated explicitly, it is feasible to numerically study ground state degeneracy and entanglement properties.

(ii) \emph{Gapless models and boundary effects.} While our discussion has primarily focused on gapped lattice models, introducing a gapless physical boundary renders the model gapless. For anyonic chain-based models, this scenario has been considered in Ref.~\cite{bottini2024haagerupsymmetry}, and similar results are expected to hold here. Additionally, the conformal field theory (CFT) perspective—including the central charge, entanglement properties, and related aspects—remains to be studied. We believe these features can be analyzed, likely via numerical methods, in analogy with Refs.~\cite{Huang2022Haagerup, corcoran2024haagerupspin, bottini2024haagerupsymmetry, Liu2023Haagerup}.

(iii) \emph{Applications in measurement-based quantum computation.} The potential use of this model for measurement-based quantum computation is another promising avenue. For finite groups, this has been partially addressed in Ref.~\cite{fechisin2023noninvertible}, but for general non-Abelian groups, Hopf algebras, and weak Hopf algebras, the problem remains largely open.

We leave these questions to future work.

\black

\vspace{1em}
\emph{Acknowledgments.} -- 
I would like to thank Dagomir Kaszlikowski for his support.
This work is supported by the National Research Foundation in Singapore, A*STAR under its CQT Bridging Grant and CQT- Return of PIs EOM YR1-10 Funding,
 and  CQT Young Researcher Career Development Grant.

\appendix
\section{Details of Haagerup boundary tube algebra $H_3$}
In this section, we provide a more detailed discussion of the Haagerup boundary tube algebra $H_3$ used in the construction of the lattice model. The first step is to enumerate all basis elements of $H_3$, which can be accomplished as follows. We group the basis elements according to the bulk strings. For the bulk string $\mathbf{1}$, there are $6 \times 6 = 36$ basis elements:
\begin{align*}
         \begin{aligned}
    \begin{tikzpicture}[scale=0.5]
     \begin{scope}
            \fill[gray!20]
                (0,1.5) arc[start angle=90, end angle=270, radius=1.5] -- 
                (0,-0.5) arc[start angle=270, end angle=90, radius=0.5] -- cycle;
        \end{scope}
         \draw[line width=.7pt,black] (0,0.5)--(0,1.5);
         \draw[line width=.7pt,black] (0,-0.5)--(0,-0.8);
         \draw[line width=.7pt,black] (0,-0.8)--(0,-1.2);
         \draw[line width=.7pt,black] (0,-1.2)--(0,-1.5);
         \draw[dotted,line width=.7pt,red] (0,0.8) arc[start angle=90, end angle=270, radius=0.8];
         \node[ line width=0.6pt, dashed, draw opacity=0.5] (a) at (-1,0){$\scriptstyle \mathbf{1}$};
        \node[ line width=0.6pt, dashed, draw opacity=0.5] (a) at (-0.3,1.3){$\scriptstyle a$};
        \node[ line width=0.6pt, dashed, draw opacity=0.5] (a) at (-0.3,-1.3){$\scriptstyle b$};
        \end{tikzpicture}
    \end{aligned},\quad \forall a,b\in \Irr(\mathcal{H}_3).
\end{align*}
Similarly, for a bulk string $\alpha \in \Irr(\cH_3)$, we need to count the number of independent bottom and top vertices. The total dimension is also $6 \times 6 = 36$ (see the fusion table, Table~\ref{tab:H3fusion}):
\begin{align*}
        \begin{aligned}
    \begin{tikzpicture}[scale=0.5]
     \begin{scope}
            \fill[gray!20]
                (0,1.5) arc[start angle=90, end angle=270, radius=1.5] -- 
                (0,-0.5) arc[start angle=270, end angle=90, radius=0.5] -- cycle;
        \end{scope}
         \draw[line width=.7pt,black] (0,0.5)--(0,0.8);
         \draw[line width=.7pt,black] (0,0.8)--(0,1.2);
         \draw[line width=.7pt,black] (0,1.2)--(0,1.5);
         \draw[line width=.7pt,black] (0,-0.5)--(0,-0.8);
         \draw[line width=.7pt,black] (0,-0.8)--(0,-1.2);
         \draw[line width=.7pt,black] (0,-1.2)--(0,-1.5);
         \draw[line width=.7pt,red] (0,0.8) arc[start angle=90, end angle=270, radius=0.8];
        \node[ line width=0.6pt, dashed, draw opacity=0.5] (a) at (-1,0){$\scriptstyle \alpha$};
        \node[ line width=0.6pt, dashed, draw opacity=0.5] (a) at (-0.3,1.3){$\scriptstyle d$};
        \node[ line width=0.6pt, dashed, draw opacity=0.5] (a) at (-0.3,0.3){$\scriptstyle c$}; 
        \node[ line width=0.6pt, dashed, draw opacity=0.5] (a) at (-0.3,-0.3){$\scriptstyle b$}; 
        \node[ line width=0.6pt, dashed, draw opacity=0.5] (a) at (-0.3,-1.3){$\scriptstyle a$};
       \node[ line width=0.6pt, dashed, draw opacity=0.5] (a) at (0.3,-1){$\scriptstyle \mu$};
        \node[ line width=0.6pt, dashed, draw opacity=0.5] (a) at (0.3,1){$\scriptstyle \nu$};
        \end{tikzpicture}
    \end{aligned} \;\;   a,b,c,d\in \Irr(\cH_3), \mu,\nu \in \Hom_{\cH_3}.
\end{align*}
By repeating this procedure for all bulk strings, the corresponding dimensions are given by
\begin{equation*}
    \begin{tabular}{c|c}
    \hline
        \hline
       bulk string  & dimension \\
         \hline
      $\mathbf{1}$   & 36 \\
      \hline
       $\alpha$   & 36 \\
      \hline
        $\alpha^2$   & 36 \\
      \hline
           $\rho$   & 225 \\
      \hline
           ${}_{\alpha}\rho$   & 225 \\
      \hline
           ${}_{\alpha^2}\rho$   & 225 \\
      \hline
          \hline
    \end{tabular}
\end{equation*}
The total dimension of $H_3$ is $783$, which is quite large; every local space in the Haagerup cluster state possesses this dimension. Consequently, numerical checks of the model's properties are computationally demanding.

If we denote these basis elements as $e_i$, the unit axiom is simply $1 \cdot e_i = e_i$ for all $i$. Multiplication is characterized by $e_i \cdot e_j =\sum_k A_{ij}^k e_k$, comultiplication by $\Delta(e_i) = \sum_{j,k} C_{i}^{jk} e_j \otimes e_k$, the counit by $\varepsilon(e_i)$, and the antipode by $S(e_i) = \sum_j S_{i}^j e_j$.
Notice that Eqs.~\eqref{eq:unitH3}--\eqref{eq:antipodeTH3} already provide all the structure constants for unit, counit, comultiplication and the antipode; we only need an explicit expression for the multiplication. 

Since $\cH_3$ is a multiplicity-free fusion category, meaning that $\operatorname{dim}\operatorname{Hom}(a \otimes b, c) = 0$ or $1$, so the vertex labels can be omitted.
To compute the structure constants for multiplication, we need to evaluate the right-hand side of Eq.~\eqref{eq:multiConst}.

\begin{align*}
  &  \begin{aligned}
    \begin{tikzpicture}
        \begin{scope}
            \fill[gray!20]
                (0,1.1) arc[start angle=90, end angle=270, radius=1.1] -- 
                (0,-0.5) arc[start angle=270, end angle=90, radius=0.5] -- cycle;
        \end{scope}
           \draw[line width=0.6pt,black,->] (0,0.5)--(0,0.7);
           \draw[line width=0.6pt,black,->] (0,0.7)--(0,1.0);
           \draw[line width=0.6pt,black] (0,0.5)--(0,1.1);
        \draw[line width=.6pt,black] (0,-0.5)--(0,-1.1);
        \draw[line width=0.6pt,black,->] (0,-1.1)--(0,-0.9);
        \draw[line width=0.6pt,black,->] (0,-1.1)--(0,-0.6);
        \draw[red, line width=0.6pt] (0,0.8) arc[start angle=90, end angle=270, radius=0.8];
           \draw[red, line width=0.6pt, ->] (0,-0.8) arc[start angle=270, end angle=180, radius=0.8];
        \node[line width=0.6pt, dashed, draw opacity=0.5] at (0,1.3) {$g$};
        \node[line width=0.6pt, dashed, draw opacity=0.5] at (0,-1.3) {$c$};
        \node[line width=0.6pt, dashed, draw opacity=0.5] at (-1,0) {$a$};
        \node[line width=0.6pt, dashed, draw opacity=0.5] at (0,-0.2) {$e$};
        \node[line width=0.6pt, dashed, draw opacity=0.5] at (0,0.2) {$f$};
    \end{tikzpicture}
\end{aligned}
\cdot
\begin{aligned}
    \begin{tikzpicture}
        \begin{scope}
            \fill[gray!20]
                (0,1.1) arc[start angle=90, end angle=270, radius=1.1] -- 
                (0,-0.5) arc[start angle=270, end angle=90, radius=0.5] -- cycle;
        \end{scope}
           \draw[line width=0.6pt,black,->] (0,0.5)--(0,0.7);
           \draw[line width=0.6pt,black,->] (0,0.7)--(0,1.0);
           \draw[line width=0.6pt,black] (0,0.5)--(0,1.1);
        \draw[line width=.6pt,black] (0,-0.5)--(0,-1.1);
        \draw[line width=0.6pt,black,->] (0,-1.1)--(0,-0.9);
        \draw[line width=0.6pt,black,->] (0,-1.1)--(0,-0.6);
        \draw[red, line width=0.6pt] (0,0.8) arc[start angle=90, end angle=270, radius=0.8];
           \draw[red, line width=0.6pt, ->] (0,-0.8) arc[start angle=270, end angle=180, radius=0.8];
        \node[line width=0.6pt, dashed, draw opacity=0.5] at (0,1.3) {$g'$};
        \node[line width=0.6pt, dashed, draw opacity=0.5] at (0,-1.3) {$c'$};
        \node[line width=0.6pt, dashed, draw opacity=0.5] at (-1,0) {$a'$};
        \node[line width=0.6pt, dashed, draw opacity=0.5] at (0,-0.2) {$e'$};
        \node[line width=0.6pt, dashed, draw opacity=0.5] at (0,0.2) {$f'$};
    \end{tikzpicture}
\end{aligned}=\delta_{e,c'}\delta_{f,g'}
\begin{aligned}
    \begin{tikzpicture}
        \begin{scope}
            \fill[gray!20]
                (0,1.7) arc[start angle=90, end angle=270, radius=1.7] -- 
                (0,-0.5) arc[start angle=270, end angle=90, radius=0.5] -- cycle;
        \end{scope}
           \draw[line width=0.6pt,black,->] (0,0.5)--(0,0.7);
           \draw[line width=0.6pt,black,->] (0,0.7)--(0,1.0);
                      \draw[line width=0.6pt,black,->] (0,0.5)--(0,1.6);
           \draw[line width=0.6pt,black] (0,0.5)--(0,1.7);
        \draw[line width=.6pt,black] (0,-0.5)--(0,-1.7);
        \draw[line width=0.6pt,black,->] (0,-1.1)--(0,-0.9);
        \draw[line width=0.6pt,black,->] (0,-1.1)--(0,-0.6);
        \draw[line width=0.6pt,black,->] (0,-1.7)--(0,-1.4);
        \draw[red, line width=0.6pt] (0,0.8) arc[start angle=90, end angle=270, radius=0.8];
           \draw[red, line width=0.6pt, ->] (0,-0.8) arc[start angle=270, end angle=180, radius=0.8];
        \draw[red, line width=0.6pt] (0,1.3) arc[start angle=90, end angle=270, radius=1.3];
           \draw[red, line width=0.6pt, ->] (0,-1.3) arc[start angle=270, end angle=180, radius=1.3];
                            \node[line width=0.6pt, dashed, draw opacity=0.5] at (0.3,1.6) {$g$};
                                      \node[line width=0.6pt, dashed, draw opacity=0.5] at (0.3,-1.6) {$c$};
        \node[line width=0.6pt, dashed, draw opacity=0.5] at (0.3,1.1) {$g'$};
        \node[line width=0.6pt, dashed, draw opacity=0.5] at (0.3,-1.1) {$c'$};
        \node[line width=0.6pt, dashed, draw opacity=0.5] at (-1,0) {$a'$};
        \node[line width=0.6pt, dashed, draw opacity=0.5] at (0,-0.2) {$e'$};
        \node[line width=0.6pt, dashed, draw opacity=0.5] at (0,0.2) {$f'$};
                \node[line width=0.6pt, dashed, draw opacity=0.5] at (-1.5,0) {$a$};
    \end{tikzpicture}
\end{aligned}\nonumber\\
=&\delta_{e,c'}\delta_{f,g'} \sum_{g'',c''} [(F_{aa'f'}^g)^{-1}]_{g'}^{g''}[(F^{aa'e'}_c)^{-1}]_{c'}^{c''}
\begin{aligned}
    \begin{tikzpicture}
        \begin{scope}
            \fill[gray!20]
                (0,1.7) arc[start angle=90, end angle=270, radius=1.7] -- 
                (0,-0.5) arc[start angle=270, end angle=90, radius=0.5] -- cycle;
        \end{scope}
           \draw[line width=0.6pt,black,->] (0,0.5)--(0,0.7);
          \draw[line width=0.6pt,black,->] (0,0.5)--(0,1.6);
           \draw[line width=0.6pt,black] (0,0.5)--(0,1.7);
        \draw[line width=.6pt,black] (0,-0.5)--(0,-1.7);
        \draw[line width=0.6pt,black,->] (0,-1.1)--(0,-0.6);
        \draw[line width=0.6pt,black,->] (0,-1.7)--(0,-1.4);
        \draw[red, line width=0.6pt] (0,1) arc[start angle=90, end angle=150, radius=1];
        \draw[red, line width=0.6pt, ->] (150:1) arc[start angle=-210, end angle=-260, radius=1];
             \draw[red, line width=0.6pt] (0,-1) arc[start angle=270, end angle=210, radius=1];
           \draw[red, line width=0.6pt, ->] (0,-1) arc[start angle=270, end angle=225, radius=1];
            \draw[red, line width=0.6pt] (210:1) arc[start angle=-90, end angle=90, radius=0.5];
             \draw[red, line width=0.6pt, ->] (210:1) arc[start angle=-90, end angle=0, radius=0.5];
             \draw[red, line width=0.6pt] (150:1) arc[start angle=90, end angle=270, radius=0.5];
             \draw[red, line width=0.6pt, ->] (210:1) arc[start angle=-90, end angle=-180, radius=0.5];
                            \node[line width=0.6pt, dashed, draw opacity=0.5] at (0.3,1.6) {$g$};
                                      \node[line width=0.6pt, dashed, draw opacity=0.5] at (0.3,-1.6) {$c$};
        \node[line width=0.6pt, dashed, draw opacity=0.5] at (-0.3,1.2) {$g''$};
        \node[line width=0.6pt, dashed, draw opacity=0.5] at (-0.3,-1.2) {$c''$};
        \node[line width=0.6pt, dashed, draw opacity=0.5] at (-0.6,0) {$a'$};
        \node[line width=0.6pt, dashed, draw opacity=0.5] at (0,-0.2) {$e'$};
        \node[line width=0.6pt, dashed, draw opacity=0.5] at (0,0.2) {$f'$};
                \node[line width=0.6pt, dashed, draw opacity=0.5] at (-1.6,0) {$a$};
    \end{tikzpicture}
\end{aligned}\nonumber\\
=&\delta_{e,c'}\delta_{f,g'} \sum_{g'',c''}[(F_{aa'f'}^g)^{-1}]_{g'}^{g''}[(F^{aa'e'}_c)^{-1}]_{c'}^{c''} \delta_{g'',c''} \sqrt{\frac{d_{a}d_{a'}}{d_{c''}}}
\begin{aligned}
    \begin{tikzpicture}
        \begin{scope}
            \fill[gray!20]
                (0,1.1) arc[start angle=90, end angle=270, radius=1.1] -- 
                (0,-0.5) arc[start angle=270, end angle=90, radius=0.5] -- cycle;
        \end{scope}
           \draw[line width=0.6pt,black,->] (0,0.5)--(0,0.7);
           \draw[line width=0.6pt,black,->] (0,0.7)--(0,1.0);
           \draw[line width=0.6pt,black] (0,0.5)--(0,1.1);
        \draw[line width=.6pt,black] (0,-0.5)--(0,-1.1);
        \draw[line width=0.6pt,black,->] (0,-1.1)--(0,-0.9);
        \draw[line width=0.6pt,black,->] (0,-1.1)--(0,-0.6);
        \draw[red, line width=0.6pt] (0,0.8) arc[start angle=90, end angle=270, radius=0.8];
           \draw[red, line width=0.6pt, ->] (0,-0.8) arc[start angle=270, end angle=180, radius=0.8];
        \node[line width=0.6pt, dashed, draw opacity=0.5] at (0,1.3) {$g$};
        \node[line width=0.6pt, dashed, draw opacity=0.5] at (0,-1.3) {$c$};
        \node[line width=0.6pt, dashed, draw opacity=0.5] at (-1,0) {$c''$};
        \node[line width=0.6pt, dashed, draw opacity=0.5] at (0,-0.2) {$e'$};
        \node[line width=0.6pt, dashed, draw opacity=0.5] at (0,0.2) {$f'$};
    \end{tikzpicture}
\end{aligned}\nonumber\\
=&\delta_{e,c'}\delta_{f,g'} \sum_{c''}[(F_{aa'f'}^g)^{-1}]_{g'}^{c''}[(F^{aa'e'}_c)^{-1}]_{c'}^{c''}  \sqrt{\frac{d_{a}d_{a'}}{d_{c''}}}
\begin{aligned}
    \begin{tikzpicture}
        \begin{scope}
            \fill[gray!20]
                (0,1.1) arc[start angle=90, end angle=270, radius=1.1] -- 
                (0,-0.5) arc[start angle=270, end angle=90, radius=0.5] -- cycle;
        \end{scope}
           \draw[line width=0.6pt,black,->] (0,0.5)--(0,0.7);
           \draw[line width=0.6pt,black,->] (0,0.7)--(0,1.0);
           \draw[line width=0.6pt,black] (0,0.5)--(0,1.1);
        \draw[line width=.6pt,black] (0,-0.5)--(0,-1.1);
        \draw[line width=0.6pt,black,->] (0,-1.1)--(0,-0.9);
        \draw[line width=0.6pt,black,->] (0,-1.1)--(0,-0.6);
        \draw[red, line width=0.6pt] (0,0.8) arc[start angle=90, end angle=270, radius=0.8];
           \draw[red, line width=0.6pt, ->] (0,-0.8) arc[start angle=270, end angle=180, radius=0.8];
        \node[line width=0.6pt, dashed, draw opacity=0.5] at (0,1.3) {$g$};
        \node[line width=0.6pt, dashed, draw opacity=0.5] at (0,-1.3) {$c$};
        \node[line width=0.6pt, dashed, draw opacity=0.5] at (-1,0) {$c''$};
        \node[line width=0.6pt, dashed, draw opacity=0.5] at (0,-0.2) {$e'$};
        \node[line width=0.6pt, dashed, draw opacity=0.5] at (0,0.2) {$f'$};
    \end{tikzpicture}
\end{aligned}
\end{align*}
The $F$-symbol is explicitly given in Ref.~\cite{osborne2019fsymbolsh3fusioncategory}; with this, all structure constants for $H_3$ can be obtained.

\black

\bibliographystyle{apsrev4-1-title}
\bibliography{Jiabib}

\end{document}